\documentclass[12pt,a4paper]{article}
\usepackage{amssymb,amsmath,color,graphicx}
\catcode`!=11
\let\!int\int \def\int{\displaystyle\!int}
\let\!lim\lim \def\lim{\displaystyle\!lim}
\let\!cap\cap \def\cap{\displaystyle\!cap}
\let\!cup\cup \def\cup{\displaystyle\!cup}
\let\!sup\sup \def\sup{\displaystyle\!sup}
\let\!inf\inf \def\inf{\displaystyle\!inf}
\let\!sum\sum \def\sum{\displaystyle\!sum}
\let\!max\max \def\max{\displaystyle\!max}
\let\!min\min \def\min{\displaystyle\!min}
\let\!frac\frac \def\frac{\displaystyle\!frac}
\catcode`!=12
\let\oldsection\section
\renewcommand\section{\setcounter{equation}{0}\oldsection}

\newtheorem{definition}{Definition}[section]
\newtheorem{lemma}{Lemma}[section]
\newtheorem{theorem}{Theorem}[section]
\newtheorem{remark}{Remark}[section]
\newtheorem{proposition}{Proposition}[section]

\newtheorem{notation}{Notation}[section]

\newcommand{\be}{\begin{equation}}
\newcommand{\ee}{\end{equation}}
\newcommand{\baa}{\begin{array}}
\newcommand{\eaa}{\end{array}}
\newcommand{\ba}{\begin{eqnarray}}
\newcommand{\ea}{\end{eqnarray}}
\allowdisplaybreaks

\begin{document}

\title{\bf Two-Dimensional Curved Fronts in a Periodic Shear Flow}

\author{Mohammad El Smaily$^{\hbox{\small{ a}}}$, Fran\c cois Hamel$^{\hbox{\small{ b}}}$ and Rui Huang$^{\hbox{\small{ c}}}$\thanks{The first author is partially supported by a PIMS postdoctoral fellowship and by an NSERC grant under the supervision of Professor Nassif Ghoussoub. The second author is partially supported by the French ``Agence Nationale de la Recherche" within the projects ColonSGS and PREFERED. He is also indebted to the Alexander von~Humboldt Foundation for its support.}\\
\\
\footnotesize{$^{\hbox{a}}$ Department of Mathematics, University of British Columbia}\\
\footnotesize{$\&$ Pacific Institute for the Mathematical Sciences}\\
\footnotesize{1984 Mathematics Road, V6T 1Z2, Vancouver, BC, Canada}\\
\footnotesize{$^{\hbox{b }}$Aix-Marseille Universit\'e \& Institut Universitaire de France}\\
\footnotesize{LATP, FST, Avenue Escadrille Normandie-Niemen, F-13397 Marseille Cedex 20, France}\\
\footnotesize{$^{\hbox{c }}$ School of Mathematical Sciences, South China Normal University,}\\
\footnotesize{Guangzhou, Guangdong 510631, China}}

\date{}

\maketitle

\begin{abstract}
\vskip 0.2cm
\noindent{}This paper is devoted to the study of travelling fronts of reaction-diffusion equations with periodic advection in the whole plane $\mathbb R^2$. We are interested in curved fronts satisfying some ``conical" conditions at infinity. We prove that there is a minimal speed $c^*$ such that curved fronts with speed~$c$ exist if and only if $c\geq c^*$. Moreover, we show that such curved fronts are decreasing in the direction of propagation, that is they are increasing in time. We also give some results about the asymptotic behaviors of the speed with respect to the advection, diffusion and reaction coefficients.\hfill\break

\noindent{\it Keywords:} \rm curved fronts, reaction-advection-diffusion equation, minimal speed, monotonicity of curved fronts.\hfill\break

\end{abstract}


\section{Introduction and main results}

In this paper, we consider the following reaction-advection-diffusion equation
\begin{eqnarray}\label{u equation}\begin{array}{l}
\displaystyle\frac{\partial u}{\partial t} =\Delta u+q(x)\frac{\partial u}{\partial y} +f(u),\;\hbox{for all} \; t\,\in\,\mathbb{R},\;(x,y)\,\in\,\mathbb{R}^2,\end{array}
\end{eqnarray}
where the advection coefficient $q(x)$ belongs to $C^{0,\delta}(\mathbb R)$ for some $\delta>0$, and satisfies
\begin{equation}\label{cq}
\forall\,x\in\mathbb{R},\quad q(x+L)=q(x) \quad \hbox{and} \quad \displaystyle{\int_{0}^L q(x)\;dx=0}
\end{equation}
for some $L>0$. The second condition for $q$ is a normalization condition. The nonlinearity$f$ is assumed to satisfy the following conditions
\begin{eqnarray}\label{cf}\left\{\begin{array}{ll}
f\; \hbox{is defined on $\mathbb{R}$, Lipschitz continuous, and}\; f\equiv0\;\hbox{in}\;\mathbb{R}\setminus{(0,1)},\vspace{3pt}\\
f\; \hbox{is a concave function of class}\; C^{1,\delta}\; \hbox{in}\; $[0,1]$,\vspace{3pt}\\
f'(0)>0,\; f'(1)<0,\hbox{ and } f(s)>0\;\hbox{ for all}\;s\in(0,1),\end{array}\right.
\end{eqnarray}
where $f$ is assumed to be right and left differentiable at~$0$ and~$1$, respectively ($f'(0)$ and~$f'(1)$ then stand for the right and left derivatives at~$0$ and~$1$). A typical example of such a function $f$ is the quadratic nonlinearity $f(u)=u(1-u)$ which was initially considered by Fisher \cite{c3fi} and Kolmogorov, Petrovsky and Piskunov \cite{c3kpp}. The equation~(\ref{u equation}) arises in various combustion and biological models, such as population dynamics and gene developments where $u$ stands for the relative concentration of some substance (see Aronson and Weinberger \cite{c3aw}, Fife \cite{c3f} and Murray \cite {c3m} for details). In combustion, equation~(\ref{u equation}) arises in models of flames in a periodic shear flow, like in simplified Bunsen flames models with a perforated burner, and~$u$ stands for the normalized temperature.

We are interested in the travelling front solutions of (\ref{u equation}) which have the form
$$u(t,x,y)=\phi(x,y+ct)$$
for all $(t,x,y)\in\mathbb R\times\mathbb R^2$, and for some positive constant $c$ which denotes the speed of propagation in the vertical direction~$-y$. Thus, we are led to the following elliptic equation
\begin{equation}\label{c equation}
\Delta\phi+(q(x)-c)\partial_{y}\phi+f(\phi)=0 \;\hbox{ for all }\;(x,y)\in\mathbb{R}^2,
\end{equation}
where the notation $\partial_{y}\phi$ means the partial derivative of the function $\phi$ with respect to the variable $y$.

We assume that the solutions $\phi$ of the equation (\ref{c equation}) are normalized so that $0\le\phi\le1$. We look in this paper for solutions of (\ref{c equation}) which satisfy the following ``conical" conditions at infinity
\begin{eqnarray}\label{cc}\left\{\begin{array}{l}
\displaystyle{\lim_{l\rightarrow-\infty}}\Big(\displaystyle{\sup_{(x,y)\in C^{-}_{\alpha,\beta,l}}}\phi(x,y)\Big)= 0,\vspace{3pt}\\
\displaystyle{\lim_{l\rightarrow+\infty}\Big(\inf_{(x,y)\in C^{+}_{\alpha,\beta,l}}\phi(x,y)\Big)= 1,}\end{array}\right.
\end{eqnarray}
where $\alpha$ and $\beta$ are given in $(0,\pi)$ such that $\alpha+\beta\le \pi$ and the lower and upper cones~$C^{-}_{\alpha,\beta,l}$ and~$C^{+}_{\alpha,\beta,l}$ are defined as follows:

\begin{definition}\label{cones}
For each real number $l$, the lower cone $C^{-}_{\alpha,\beta,l}$ is defined by
$$\begin{array}{cl}
&C^{-}_{\alpha,\beta,l}=\big\{(x,y)\in\mathbb{R}^2,~~ y\leq x\cot\alpha+l~\hbox{ whenever }~ x\leq0\vspace{3pt}\\
&\hskip3cm\hbox{ and }~ y\leq-x\cot\beta+l ~\hbox{ whenever }~ x\geq0\big\}\end{array}$$
and then the upper cone $C^{+}_{\alpha,\beta,l}$ is defined by
$$C^{+}_{\alpha,\beta,l}=\overline{\mathbb{R}^2 \setminus C^{-}_{\alpha,\beta,l}},$$
see Figure $\ref{conesf}$ for a geometrical description.
\begin{figure}
\begin{center}
\includegraphics[width=3.5in]{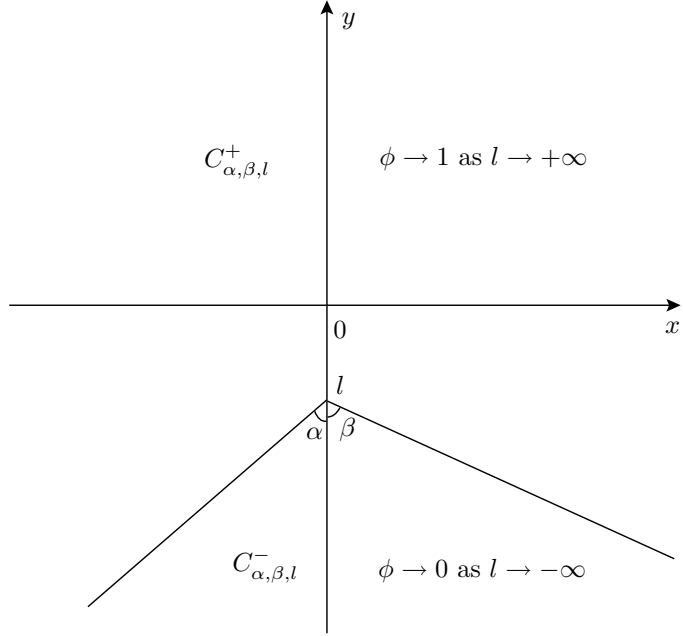}\\
\caption{the lower and upper cones $C^{-}_{\alpha,\beta,l}$ and $C^{+}_{\alpha,\beta,l}$.}\label{figc3}
\label{conesf}\end{center}
\end{figure}
\end{definition}

Because of the strong elliptic maximum principle, a solution $\phi$ of the equation (\ref{c equation}) that is defined in the whole plane $\mathbb R^2$ and satisfies $0\le\phi\le1$, is either identically equal to~$0$ or~$1$, or $0<\phi(x,y)<1$ for all $(x,y)\in\mathbb R^2$. By the ``conical" conditions at infinity~(\ref{cc}), only the case of $0<\phi(x,y)<1$ for all $(x,y)\in\mathbb R^2$ will then be considered in the present paper.

In order to motivate our study, let us first recall a very simple case of travelling fronts for the reaction-diffusion (with no advection) equation
\begin{equation}\label{nonadvection}
\frac{\partial u}{\partial t}-\Delta u=f(u)
\end{equation}
in the whole plane $\mathbb R^2$. It is well known from \cite{c3kpp} that for any $c\geq2\sqrt{f'(0)},$ the above equation has a planar travelling front moving in an arbitrarily given unit direction $-e$, having the form  $u(t,x)=\phi(x\cdot e+ct)$ and satisfying the conditions $\phi(-\infty)=0$ and $\phi(+\infty)=1$. Recently, the problems about curved travelling fronts of the reaction-diffusion (with no advection) equation (\ref{nonadvection}) equipped with the conical conditions at infinity of type~(\ref{cc}) with $\alpha=\beta$ have been the subject of intensive study by many authors, for various types of nonlinearities. For example, Bonnet and Hamel \cite{c3bf99} considered such type of problems with a ``combustion'' nonlinearity $f$, namely,
$$\exists \,\theta\in(0,1),\ f=0\ \hbox{on}\ [0,\theta]\ \hbox{and}\ f'(1)<0,$$
which comes from the model of premixed bunsen flames. They proved the existence of curved travelling fronts and gave an explicit formula that relates the speed of propagation and the angle of the tip of the flame. One can also find some generalizations of the above results and further qualitative properties in \cite{c3hm00,c3hmr04}. For the case of bistable nonlinearity~$f$ satisfying
\begin{equation*}\begin{cases}
\exists\,\theta\in(0,1),\ f(0)=f(\theta)=f(1)=0, \ f'(0)<0,\ f'(1)<0,\ f'(\theta)>0,\\
f<0\  \hbox{on}\  (0,\theta)\cup(1,+\infty),\ f>0 \ \hbox{on}\ (-\infty,0)\cup(\theta,1),\end{cases}
\end{equation*}
Hamel, Monneau and Roquejoffre \cite{c3hmr05,c3hmr06} and Ninomiya and Taniguchi~\cite{c3nt05,c3nt06} proved existence and uniqueness results and qualitative properties of such kind of conical fronts (see also \cite{rrm,t1,t2} for further stability results and the study of pyramidal fronts). For KPP nonlinearities, conical and more general curved fronts are also known to exist for equation~(\ref{nonadvection}) (see \cite{c3hn01}).

In addition to the above mentioned literature, some works have been devoted to the study of the reaction-advection-diffusion equations of the type (\ref{u equation}). A well-known paper about this issue is the one by Berestycki and Nirenberg \cite{c3bn92}, where the authors set the reaction-advection-diffusion equation in a straight infinite cylinder and consider the travelling fronts of the reaction-advection-diffusion equation satisfying Neumann no-flux conditions on the boundary of the cylinder and approaching $0$ and $1$ at both infinite sides of the cylinder respectively. Later, Berestycki and Hamel \cite{c3bh02} and Weinberger~\cite{c3wein02} investigated reaction-diffusion equations with periodic advection in a very general framework, and proved the existence of pulsating travelling fronts (some of their results will be recalled below).

However, as far as we know, except recent works of Haragus and Scheel \cite{c3hs1,c3hs2} on some equations of the type (\ref{c equation}) with $\alpha$ and $\beta$ close to $\pi/2$, the reaction-advection-diffusion equation of type (\ref{u equation}) and its corresponding elliptic equation (\ref{c equation}) equipped with conical conditions (\ref{cc}) have not been studied yet for general angles $\alpha$ and $\beta$ and for general periodic shear flow. The purpose of this paper is to prove the existence, nonexistence and monotonicity results for the solutions of the semilinear elliptic equation (\ref{c equation}) with the non-standard conical conditions at infinity (\ref{cc}). In fact, the main difficulties in the present paper arise from these conical conditions at infinity and from the fact that the domain is not compact in the direction orthogonal to the direction of propagation.

Before stating our main results of this paper, we first give the following notations.

\begin{notation}\label{not 1}{\rm
Let $\gamma\in(0,{\pi}/{2})$, $q=q(X)$ and $f=f(u)$ be two functions satisfying $(\ref{cq})$ and $(\ref{cf})$ respectively. Let $M=(m_{ij})_{1\le i,j\le2}$ be a positive definite symmetric matrix, that is
\be\label{Mdef+}
\exists\, c_1>0, \; \forall\,\xi\in\mathbb R^2,\; \sum_{1\le i,j\le 2}m_{i,j}\xi_i\xi_j\ge c_1|\xi|^2,
\ee
where $|\xi|^2=\xi_1^2+\xi_2^2$ for any $\xi=(\xi_1,\xi_2)\in\mathbb R^2$. Throughout this paper, $c^*_{M,q\sin\gamma,f}>0$ denotes the minimal speed of propagation of travelling fronts $0\le u\le1$ in the direction~$-Y$ in the variables $(X,Y)$ for the following reaction-advection-diffusion problem
\begin{equation}\label{subbht}\begin{array}{rcl}
\frac{\partial u}{\partial t} & \!\!=\!\! & {\rm div}(M \nabla u)+q(X)\sin\gamma \frac{\partial u}{\partial Y}\!+\!f(u),\ \, t\in \mathbb R,\; (X,Y)\!\in\mathbb R^2,\vspace{3pt}\\
u(t\!+\!\tau,X\!+\!L,Y) & \!\!=\!\! & u(t\!+\!\tau,X,Y)=u(t,X,Y\!+\!c\tau),\ \,(t,\tau,X,Y)\!\in\mathbb R^2\!\times\!\mathbb R^2,\vspace{3pt}\\
u(t,X,Y) & \!\!\!\!\underset{Y\rightarrow -\infty}{\longrightarrow}\!\!\!\! &  0,\quad u(t,X,Y)\underset{Y\rightarrow +\infty}{\longrightarrow} 1,\end{array}
\end{equation}
where the above limits hold locally in $t$ and uniformly in $X$. In other words, such fronts exist if and only if $c\ge c^*_{M,q\sin\gamma,f}$. The existence of this miminal speed $c^*_{M,q\sin\gamma,f}$ and further qualitative properties of such fronts, for even more general periodic equations, follow from {\rm{\cite{c3bh02,c3h08,c3hrpreprint,c3wein02}}} $($see also {\rm{\cite{c3bn92}}} for problems set in infinite cylinders$)$.}
\end{notation}

Our first main result in this paper is the following

\begin{theorem}\label{c3theorem1}
Let $q(x)$ be a globally $C^{0,\delta}(\mathbb R)$ function $($for some $\delta>0)$ satisfying $(\ref{cq})$. Let $f$ be a nonlinearity fulfilling $(\ref{cf})$. Then, for any given $\alpha$ and $\beta$ in $(0,\pi)$ such that $\alpha+\beta\le\pi$, there exists a positive real number $c^*$ such that\par
i$)$ for each $c\ge c^*$, the problem $(\ref{c equation})$-$(\ref{cc})$ admits a solution $(c,\phi)$;\par
ii$)$ if $c<c^*$, the problem $(\ref{c equation})$-$(\ref{cc})$ has no solution $(c,\phi)$.\\
Moreover, under the Notation $\ref{not 1}$, the value of $c^*$ is given by
\begin{equation}\label{c^*}
c^*=\max \left(\frac{c^*_{A,q\sin\alpha,f}}{\sin\alpha},\frac{c^*_{B,q\sin\beta,f}}{\sin\beta}\right),
\end{equation}
where
\be\label{defAB}
A=\left[\begin{array}{cc}1& -\cos\alpha\\ -\cos\alpha&1\end{array}\right]\quad \hbox{and}\quad B=\left[\begin{array}{cc} 1&\cos\beta\\ \cos\beta&1\end{array}\right].
\ee
\end{theorem}

Our second main result is concerned with the monotonicity of the fronts in the direction of propagation.

\begin{theorem}\label{phi is y increasing}
Under the assumptions of Theorem $\ref{c3theorem1}$, if a pair $(c,\phi)$ solves the pro\-blem~$(\ref{c equation})$-$(\ref{cc})$, then $\partial_y\phi(x,y)>0$ for all $(x,y)\in\mathbb R^2$. Consequently, the travelling front solution $u(t,x,y)=\phi(x,y+ct)$ of $(\ref{u equation})$ is increasing in time $t$.
\end{theorem}

\begin{remark}{\rm For the case of $\alpha=\beta=\pi/2$, the above results have been proved in \cite{c3bh02,c3bn92,c3wein02}, in which case $c^*=c^*_{I,q,f}$ is the minimal speed of travelling fronts for problem (\ref{subbht}) with identity matrix $M=I$. The interest of the present work is to generalize them to the case of conical asymptotic conditions (\ref{cc}) with angles $\alpha$ and $\beta$ which may be smaller or larger than $\pi/2$. The condition $\alpha+\beta\le\pi$ which is used in the construction of the fronts can be viewed as a global concavity of the level sets of the fronts with respect to the variable $y$. It is unclear that this condition is necessary in general. Actually, it follows from Section 4 that Theorem \ref{phi is y increasing} still holds for any $\alpha$ and $\beta$ in $(0,\pi)$.}
\end{remark}

The value of $c^*$ in Theorem \ref{c3theorem1} is given in terms of the known minimal speeds of propagation of ``planar'' pulsating travelling fronts for two auxiliary (left and right) problems of type (\ref{subbht}). Throughout the paper, we use the word ``planar" to mean that, for pro\-blem~(\ref{subbht}), any level set of~$u$ is trapped between two parallel planes. A rigorous result about the existence of the minimal speed of propagation of pulsating travelling fronts in general periodic domains was given in \cite{c3bh02}. Several variational formul{\ae} for the minimal speed of propagation have been given by Berestycki, Hamel and Nadirashvili \cite{BHN1}, El Smaily~\cite{El Smaily min max} and Weinberger~\cite{c3wein02}. Much work has been devoted to the study of the dependence of the ``planar'' minimal speed on the advection, diffusion, reaction and the geometry of the domain (see e.g. \cite{BHN1,BHN2,El Smaily,EHR,hnr,Heinze convection,Nadin,Nadin2,ryzhikzlatos,zlatos}).

In the following theorem, we study the behaviors of the \emph{conical} minimal speed $c^*$ of Theorem \ref{c3theorem1} in some asymptotic regimes and we obtain a result about the homogenized speed. To make the presentation simpler, we introduce a general notation for the conical minimal speed: given an advection $q$ and a reaction $f$ satisfying (\ref{cq}) and (\ref{cf}) respectively, and given an arbitrary $\rho>0$, we consider  the problem
\begin{equation}\label{c equation generic}
\rho\Delta\phi+(q(x)-c)\partial_{y}\phi+f(\phi)=0\;\hbox{ for all }\;(x,y)\in\mathbb{R}^2,
\end{equation}
with the conical conditions (\ref{cc}) and we denote by $c^*(\rho,q,f)$ the conical minimal speed of problem (\ref{c equation generic})-(\ref{cc}), whose existence follows from Theorem \ref{c3theorem1}. In other words, a solution~$(c,\phi)$ of (\ref{c equation generic}) satisfying (\ref{cc}) exists if and only if $c\geq c^*(\rho,q,f)$. Furthermore, it follows from Theorem \ref{c3theorem1} that  the ``conical'' minimal speed can be expressed in terms of the ``left and right'' planar minimal speeds as follows
\begin{equation}\label{generic notation c*}
c^*(\rho,q,f)=\max \left(\frac{c^*_{\rho A,q\sin\alpha,f}}{\sin\alpha},\frac{c^*_{\rho B,q\sin\beta,f}}{\sin\beta}\right).
\end{equation}
In the above notation of conical minimal speed, we use the brackets (i.e. $c^*(\cdot,\cdot,\cdot)$) while subscripts are used in the notation of the ``planar" minimal speed.

\begin{theorem}\label{c3 asymptotics}
Let $\alpha$ and $\beta$ be in $(0,\pi)$ such that $\alpha+\beta\le\pi$. Assume that the function~$f$ fulfills $(\ref{cf})$ and that the advection $q$ is a globally $C^{0,\delta}(\mathbb R)$ function $($for some $\delta>0)$ satisfying $(\ref{cq})$.\par
i$)$ {\rm Large diffusion or small reaction with a not too large/sufficiently small advection}. For each $\rho>0,$ we have
\begin{equation}\label{limit with small reaction}
\forall\,\gamma\geq 1/2,~~ \lim_{m\rightarrow0^+}\frac{c^*(\rho,m^\gamma q,mf)}{\sqrt{m}}=\frac{2\sqrt{\rho f'(0)}}{\min(\sin\alpha,\sin\beta)},
\end{equation}
and
\begin{equation}\label{limit with large diffusion}
\forall\,0\leq\gamma\leq 1/2,~~ \lim_{m\rightarrow+\infty}\frac{c^*(m\rho,m^\gamma q,f)}{\sqrt{m}}=\frac{2\sqrt{\rho f'(0)}}{\min(\sin\alpha,\sin\beta)}.
\end{equation}\par
ii$)$ {\rm Large advection}. For each $\rho>0,$ the following limit holds
\begin{equation}\label{within large advection}
\lim_{m\rightarrow+\infty}\frac{c^*(\rho,mq,f)}{m}=\max_{\substack{w\in H^1_{loc}(\mathbb{R})\backslash\{0\},\,L-\text{periodic},\\ \rho\|w'\|_{L^2(0,L)}^2\le f'(0)\|w\|_{L^2(0,L)}^2}}\frac{\int_0^Lq\,w^2}{\int_0^Lw^2}.
\end{equation}
Moreover,
\begin{equation}\label{large advection small reaction}\baa{rcl}
\lim_{\varepsilon\rightarrow0^+}\left(\lim_{m\rightarrow+\infty}\displaystyle{\frac{\displaystyle{c^{*}( \rho,mq,\varepsilon f)}}{m\sqrt{\varepsilon}}}\right) & = & \lim_{\mu\rightarrow+\infty}\left(\lim_{m\rightarrow+\infty}\displaystyle{\frac{\displaystyle{c^{*}(\mu \rho,mq,f)\times\sqrt{\mu}}}{m}}\right)\vspace{3pt}\\
& = & \displaystyle{\frac{2\sqrt{f'(0)}}{\sqrt{\rho L}}}\times\max_{w\in H^1_{loc}(\mathbb{R})\backslash\{0\},\,L-\text{periodic}}\ \frac{\int_0^Lq\,w}{\|w'\|_{L^2(0,L)}}\eaa
\end{equation}
and
\be\label{mqmuf}
\lim_{\varepsilon\rightarrow0^+}\left(\lim_{m\rightarrow+\infty}\displaystyle{\frac{\displaystyle{c^{*}( \varepsilon,mq,f)}}{m}}\right)=\lim_{\mu\rightarrow+\infty}\left(\lim_{m\rightarrow+\infty}\displaystyle{\frac{\displaystyle{c^{*}(\rho,mq,\mu f)}}{m}}\right)=\max_{[0,L]}\,q.
\ee\par
iii$)$ {\rm Homogenized speed}. Assume here that $q$ is $1$-periodic and its average is zero. For each $L>0$, let $q_{_L}(x)=q\left({x}/{L}\right)$ for all $x\in\mathbb{R}$. Then, for each $\rho>0$,
\begin{equation}\label{limit as L goes to zero}
\lim_{L\rightarrow0^+}c^{*}(\rho,q_L,f)=\frac{2\sqrt{\rho f'(0)}}{\min(\sin\alpha,\sin\beta)}.
\end{equation}
\end{theorem}

\noindent{\bf{Outline of the rest of the paper.}} This paper is organized as follows. In Section \ref{existence section}, we prove the existence of a curved traveling front to the problem (\ref{c equation})-(\ref{cc}) whenever the speed $c\ge c^*$ (the first part of Theorem~\ref{c3theorem1}). In Section \ref{nonexistence section}, using some results about spreading phenomena, we prove that the problem (\ref{c equation})-(\ref{cc}) has no solution $(c,\phi)$ as soon as $c< c^*$ (the second part of Theo\-rem~\ref{c3theorem1}). In Section \ref{monotonocity section}, we first establish a generalized comparison principle for some elliptic equations in unbounded domains having the form of ``upper cones''. Then, we give the proof of Theorem \ref{phi is y increasing} by using this generalized comparison principle together with suitable estimates of the quantity ${\partial_y \phi}/{\phi}$ in lower cones and with some sliding techniques on the solutions in the $y$-variable. Lastly, Section\ref{proofs of asymptotics} is concerned with the proof of Theorem~\ref{c3 asymptotics}.


\section{Existence of a curved front $(c,\phi)$ for all $c\ge c^*$}\label{existence section}

In this section, we prove the existence of a curved front $(c,\phi)$ to the problem (\ref{c equation})-(\ref{cc}) whenever $c\ge c^*$ (the first item of Theorem \ref{c3theorem1}). The main tool is the sub/super-solution method. Roughly speaking, we construct a subsolution and a supersolution for our pro\-blem by mixing, in different ways, two pulsating travelling fronts coming from opposite sides (left and right) and having different angles with respect to the vertical axis but having the same vertical speed in some sense.\hfill\break

\noindent{\bfseries Proof of part i) of Theorem \ref{c3theorem1}.} We perform this proof in two steps.\hfill\break

\noindent{\it Step 1: Construction of a subsolution.} For any given $\gamma\in(0,\pi)$, any smooth function~$q$ sa\-tisfying (\ref{cq}), any nonlinearity $f$ fulfilling (\ref{cf}) and any constant matrix $M=(m_{ij})_{1\le i,j\le2}$ satisfying~(\ref{Mdef+}), we consider the problem (\ref{subbht}). It follows, from Theorem 1.14 in  \cite{c3bh02} that there exists a minimal speed $c^*_{M, q\sin\gamma,f}$ such  that the problem~(\ref{subbht}) admits a pulsating travelling front $(c,u)$ for each $c\ge c^*_{M, q\sin\gamma,f}$ and no solution for $c< c^*_{M, q\sin\gamma,f}$. Moreover, it is known that any such front $u$ is increasing in $t$. For any solution $(c,u)$ of the pro\-blem~(\ref{subbht}), if we denote $u(t,X,Y)=\varphi(X,Y+ct)$, then the pair $(c,\varphi)$ solves the following problem
\begin{equation}\label{subbhc}\begin{cases}
\hbox{div}(M \nabla \varphi)+(q(X)\sin\gamma-c){\partial_Y\varphi}+f(\varphi)=0, \quad (X,Y)\in\mathbb R^2,\\
\varphi(X,Y)\underset{Y\rightarrow -\infty}{\longrightarrow} 0,\quad \varphi(X,Y)\underset{Y\rightarrow +\infty}{\longrightarrow} 1,\quad \hbox{uniformly in }X\in\mathbb R,\\
\varphi(X+L,Y)=\varphi(X,Y), \quad (X,Y)\in \mathbb R^2.\end{cases}
\end{equation}
Since $u$ is increasing in $t$, we conclude that $\varphi$ is increasing in its second variable, namely~$Y$.

For any given $0<\alpha,\beta<\pi$ such that $\alpha+\beta\le\pi$, we define the matrices~$A$ and~$B$ as in~(\ref{defAB}). By choosing $M=A$ and $\gamma=\alpha$ in (\ref{subbhc}), then there exists a positive constant~$c^*_{A,q\sin\alpha,f}$ such that the problem (\ref{subbhc}) admits a solution $(c_\alpha, \varphi_{\alpha})$ if and only if $c_\alpha\ge c^*_{A,q\sin\alpha,f}$. Similarly, if we choose $M=B$ and $\gamma=\beta$ in (\ref{subbhc}), then there exists a positive constant $c^*_{B,q\sin\beta,f}$ such that the problem (\ref{subbhc}) admits a solution $(c_\beta, \varphi_{\beta})$ if and only if $c_\beta\ge c^*_{B,q\sin\beta,f}$. Consequently, for a given $c\ge c^*$, where $c^*$ is defined by (\ref{c^*}), there exist $(c_\alpha, \varphi_{\alpha})$ and $(c_\beta, \varphi_{\beta})$ as above and such that
\begin{equation}\label{subc}
c=\frac{c_\alpha}{\sin\alpha}=\frac{c_\beta}{\sin\beta}\ge c^*.
\end{equation}\par
Now, we give a candidate for a subsolution of the problem (\ref{c equation})-(\ref{cc}) as follows
\begin{equation}\label{subsolution}
\underline{\phi}(x,y)=\max\left(\varphi_{\alpha}(x,-x\cos\alpha+y\sin\alpha),\varphi_{\beta}(x,x\cos\beta+y\sin\beta)\right).
\end{equation}
In fact, by (\ref{subbhc}), it is easy to verify that $(c,\underline{\phi})$ defined by (\ref{subc}) and (\ref{subsolution}) is a subsolution of the equation (\ref{c equation}). Indeed, both functions in the max solve (\ref{c equation}). For instance, if we set $\phi_1(x,y)=\varphi_\alpha(x,-x\cos\alpha+y\sin\alpha)$, then
\begin{align*}
\Delta \phi_1+(q(x)-c)\partial_y\phi_1+f(\phi_1)=\hbox{div}(A\nabla\varphi_\alpha)+(q(x)-c)\sin\alpha\partial_Y\varphi_\alpha+f(\varphi_\alpha)=0
\end{align*}
in $\mathbb R^2$, where the quantities involving $\varphi_\alpha$ are taken values at the point $(x,-x\cos\alpha+y\sin\alpha)$. Moreover, by construction and since $\alpha+\beta\le\pi$, we know that $\underline{\phi}$ satisfies the ``conical" conditions at infinity (\ref{cc}).\hfill\break

\noindent{\it Step 2: Construction of a supersolution.} As we have done in the first step, for any $c\ge c^*$, we consider the same front $(c_{\alpha},\varphi_{\alpha})$ as in step 1, which solves the problem (\ref{subbhc}) for $M=A$ and $\gamma=\alpha$, and the same front $(c_{\beta},\varphi_{\beta})$ as in step 1, which solves the problem (\ref{subbhc}) for $M=B$ and $\gamma=\beta$ such that (\ref{subc}) holds. We claim that the following function
\begin{equation}\label{supersolution}
\bar{\phi}(x,y)=\min\left(\varphi_{\alpha}(x,-x\cos\alpha+y\sin\alpha)+\varphi_{\beta}(x,x\cos\beta+y\sin\beta),1\right)
\end{equation}
is a supersolution of the equation (\ref{c equation}). Obviously, we only need to check the case of $\varphi_{\alpha}(x,-x\cos\alpha+y\sin\alpha)+\varphi_{\beta}(x,x\cos\beta+y\sin\beta)\le 1$.

We first notice that a function $f=f(s)$ that satisfies the conditions (\ref{cf}) is sub-additive in the interval $[0,1]$. That is
$$f(s+t)\le f(s)+f(t),\quad \hbox{for all}\quad 0\le s,t\le1.$$
When $\bar\phi\le1$, then by (\ref{subbhc}), we have,
$$\baa{rcl}
\Delta \bar\phi+(q(x)-c)\partial_y\bar\phi+f(\bar\phi) & = & f(\varphi_{\alpha}+\varphi_{\beta})+\hbox{div}(A\nabla\varphi_\alpha)+(q(x)-c)\sin\alpha\,\partial_Y\varphi_\alpha\vspace{3pt}\\
& & +\hbox{div}(B\nabla\varphi_\beta)+(q(x)-c)\sin\beta\,\partial_Y\varphi_\beta\vspace{3pt}\\
& = & f(\varphi_{\alpha}+\varphi_{\beta})-f(\varphi_\alpha)-f(\varphi_\beta)\\
& \le & 0,\eaa$$
where the quantities involving $\varphi_\alpha$ (resp. $\varphi_\beta$) are taken values at the point $(x,-x\cos\alpha+y\sin\alpha)$ (resp. $(x,x\cos\beta+y\sin\beta)$). Thus, $(c,\bar\phi)$ is a supersolution of the equation (\ref{c equation}). Furthermore, the function $\bar \phi$ satisfies the conical conditions~(\ref{cc}) at infinity since $\alpha+\beta\le\pi$.

Finally, since $0\le\underline\phi\le\bar\phi\le1$ in $\mathbb R^2$, we conclude that, for any $c\ge c^*,$ the problem (\ref{c equation})-(\ref{cc}) admits a curved front $(c,\phi)$ such that $\underline\phi\le\phi\le\bar\phi$. The proof of part~i) of Theorem~\ref{c3theorem1} is then complete.\hfill$\Box$\break

Notice that it follows from the above construction that $\phi$ is close to the oblique ``planar" fronts $\varphi_{\alpha}(x,-x\cos\alpha+y\sin\alpha)$ and $\varphi_{\beta}(x,x\cos\beta+y\sin\beta)$ asymptotically on the ``left" and ``right". More precisely,
$$\lim_{A\rightarrow-\infty}\Big(\sup_{y\le x\cot\alpha+A}|\phi(x,y)-\varphi_\beta(x,x\cos\beta+y\sin\beta)|\Big)=0$$
and
$$\lim_{A\rightarrow-\infty}\Big(\sup_{y\le -x\cot\beta+A}|\phi(x,y)-\varphi_\alpha(x,-x\cos\alpha+y\sin\alpha)|\Big)=0.$$

\begin{remark}{\rm To complete this section, consider here the special ``symmetric" case. Namely, under the notations of Theorem \ref{c3theorem1}, assume $\alpha=\beta$ and $q(x)=q(-x)$ for all $x\in\mathbb R$. Then we claim that
$$c^*=\frac{c^*_{A,q\sin\alpha,f}}{\sin\alpha}=\frac{c^*_{B,q\sin\beta,f}}{\sin\beta}.$$
Indeed, let $(c^*_{A,q\sin\alpha,f}, \varphi^*_{\alpha}(X,Y))$ be a solution of the following problem
\begin{equation}\label{subbhc2}\begin{cases}
\hbox{\rm div}(A \nabla \varphi^*_\alpha(X,Y))\!+\!(q(X)\sin\alpha\!-\!c^*_{A,q\sin\alpha,f}){\partial_Y \varphi^*_\alpha(X,Y)}\!+\!f(\varphi^*_\alpha(X,Y))\!=\!0\,\hbox{in}\,\mathbb R^2,\!\!\vspace{3pt}\\
\varphi^*_\alpha(X,Y)\underset{Y\rightarrow -\infty}{\longrightarrow} 0,\quad \varphi^*_\alpha(X,Y)\underset{Y\rightarrow +\infty}{\longrightarrow} 1,\quad \hbox{uniformly in }X\in\mathbb R.\end{cases}
\end{equation}
Define $\psi(X,Y):=\varphi^*_{\alpha}(-X,Y)$ for all $(X,Y)\in\mathbb R^2$. Since $\alpha=\beta$ and $q(X)=q(-X)$ for all $X\in\mathbb R,$ then the pair $(c^*_{A,q\sin\alpha,f}, \psi)$ is a solution of the following problem
\begin{equation}\label{subbhc3}\begin{cases}
\hbox{\rm div}(B \nabla \psi(X,Y))\!+\!(q(X)\sin\alpha\!-\!c^*_{A,q\sin\alpha,f})\,{\partial_Y \psi(X,Y)}\!+\!f(\psi(X,Y))=0\hbox{ in }\mathbb R^2,\vspace{3pt}\\
\psi(X,Y)\underset{Y\rightarrow -\infty}{\longrightarrow} 0,\quad \psi(X,Y)\underset{Y\rightarrow +\infty}{\longrightarrow} 1,\quad \hbox{uniformly in }X\in\mathbb R.\end{cases}
\end{equation}
It follows from \cite{c3bh02} that $c^*_{A,q\sin\alpha,f}$ is not smaller than the minimal speed of propagation corresponding to the reaction-advection-diffusion equation having $B$ as the diffusion matrix, $q\sin\alpha=q\sin\beta$ as the advection and $f$ as the reaction term. That is, $c^*_{A,q(x)\sin\alpha,f}\ge c^*_{B,q(x)\sin\beta,f}.$ Similarly, we can prove $c^*_{B,q\sin\beta,f}\ge c^*_{A,q\sin\alpha,f}$ which leads to the equality between these two minimal speeds.}
\end{remark}


\section{Nonexistence of conical fronts $(c,\phi)$ for $c< c^*$}\label{nonexistence section}

In this section, we prove that the problem (\ref{c equation})-(\ref{cc}) has no solution $(c,\phi)$ if $c< c^*$ (the second item of Theorem \ref{c3theorem1}).
The proof mainly lies on a spreading result given by Weinberger \cite{c3wein02}.\hfill\break

\noindent{\bfseries Proof of part ii) in Theorem \ref{c3theorem1}.} Suppose to the contrary that the problem (\ref{c equation})-(\ref{cc}) admits a solution $\phi$ with a speed $c<c^*$, where $c^*$ is the value defined in (\ref{c^*}). Without loss of generality, we can assume that
$$c^*=\frac{c^*_{A,q\sin\alpha,f}}{\sin\alpha}\ge\frac{c^*_{B,q\sin\beta,f}}{\sin\beta}.$$
Under this assumption, there exists a positive constant $d$ such that
\begin{equation}\label{c3s3d}
c\sin\alpha<d<c^*_{A,q\sin\alpha,f}.
\end{equation}

Write $\phi(x,y)=\varphi(x,-x\cos\alpha+y\sin\alpha)$ for all $(x,y)\in\mathbb R^2$. Then, the function $\varphi(X,Y)$ is well defined and it solves the following equation
\begin{equation}\label{c3s3e1}
\hbox{div}(A \nabla \varphi)+(q(X)-c)\sin\alpha{\partial_Y \varphi}+f(\varphi)=0,
\quad\hbox{for all }(X,Y)\in\mathbb R^2,
\end{equation}
where $A$ is the matrix defined in the second section. Moreover, it follows from the definition of $\varphi$ and the ``conical" conditions at infinity (\ref{cc}) that
\begin{equation}\label{c3s3e2}
\lim_{Y\rightarrow-\infty}\Big(\sup_{(X,Y)\in\mathbb R^2,\,{X\le0}}\varphi(X,Y)\Big)=0.
\end{equation}
We mention that taking the supremum in the above limit over the set $\{ X\leq0\}$ is just to insure that $(X,Y)$ stays in $C_{\alpha,\beta,l}^-$ for some $l$ which goes to $-\infty$ as $Y\rightarrow-\infty$ and as a consequence we can use the conical conditions. If we let $u(t,X,Y)=\varphi(X,Y+ct\sin\alpha)$, then by (\ref{c3s3e1}), the function $u$ solves the following parabolic equation
\begin{equation}\label{c3s3e3}
\frac{\partial u}{\partial t}=\hbox{div}(A\nabla u)+q(X)\sin\alpha\frac{\partial u}{\partial Y}+f(u),\hbox{ for all }(t,X,Y)\in\mathbb R\times\mathbb R^2.
\end{equation}

Let $\hat u_0(X,Y)$ be a function of class $C^{0,\mu}(\mathbb R^2)$ (for some positive $\mu$) such that
\begin{equation}\label{c3s3e4}\left\{\begin{array}{rl}
\forall\,X\in\mathbb R,~\forall\,Y\leq0,&\hat u_0(X,Y)=0,\vspace{4 pt} \\
\exists\, Y_0>0,& \inf_{(X,Y)\in\mathbb R^2,\,Y\ge Y_0}\hat u_0(X,Y)>0, \vspace{4 pt}\\
\forall\, (X,Y)\in\mathbb R^2,&0\le\hat u_0(X,Y)\le u(0,X,Y).\end{array}\right.
\end{equation}
Let $\hat u(t,X,Y)$ be a classical solution of the following Cauchy problem
$$\begin{cases}
\frac{\partial \hat u}{\partial t}=\hbox{div}(A\nabla \hat u)+q(X)\sin\alpha\frac{\partial \hat u}{\partial Y}+f(\hat u),\ \hbox{for all}\  t>0,\ (X,Y)\in\mathbb R^2,\vspace{4pt}\\
\hat u(0,X,Y)=\hat u_0(X,Y),\ \hbox{for all}\ (X,Y)\in\mathbb R^2.\end{cases}$$
Under the conditions (\ref{c3s3e4}) on $\hat u_0$ and the assumptions (\ref{cf}) on the nonlinearity $f$, the results of Weinberger \cite{c3wein02} imply that for any given $r>0$, we have
$$\lim_{t\rightarrow+\infty}\sup_{|Y|\le r,X\in\mathbb R}\hat u(t,X,Y-c't)=0,\ \hbox{for each}\ c'>c^*_{A,q\sin\alpha,f}$$
and
\begin{equation}\label{c3s3e6}
\lim_{t\rightarrow+\infty}\inf_{|Y|\le r,X\in\mathbb R}\hat u(t,X,Y-c't)=1,\ \hbox{for each}\ c'<c^*_{A,q\sin\alpha,f}.
\end{equation}

On the other hand, since $0\le\hat u(0,X,Y)\le u(0,X,Y)$ in $\mathbb R^2$ and both $u$ and $\hat u$ solve the same parabolic equation (\ref{c3s3e3}), the parabolic maximum principle implies that
\begin{equation}\label{c3s38}
\hat u(t,X,Y)\le u(t,X,Y)\quad \hbox{for all}\  t\ge0\ \hbox{and}\  (X,Y)\in\mathbb R^2.
\end{equation}
The assumption that $(c\sin\alpha-d)<0$ implies that $Y+(c\sin\alpha-d)t\rightarrow-\infty$ as $t\rightarrow+\infty$ for $|Y|\leq r$. We conclude from (\ref{c3s3d}), (\ref{c3s3e2}) and (\ref{c3s38}) that for any $r>0$, all limits below exist and
\begin{align*}
0\le\lim_{t\rightarrow+\infty}\inf_{|Y|\le r,\,X\in\mathbb R}\hat u(t,X,Y-dt)&\le\lim_{t\rightarrow+\infty}\inf_{|Y|\le r,\,X\le0}\hat u(t,X,Y-dt)\\
&\le\lim_{t\rightarrow+\infty}\inf_{|Y|\le r,\,X\le0}u(t,X,Y-dt)\\
&=\lim_{t\rightarrow+\infty}\inf_{|Y|\le r,\,X\le0}\varphi(X,Y+(c\sin\alpha-d)t)\\
&\le\lim_{t\rightarrow+\infty}\sup_{|Y|\le r,\,X\le0}\varphi(X,Y+(c\sin\alpha-d)t)\\
&=0,
\end{align*}
which contradicts (\ref{c3s3e6}) with $c'=d$ and eventually completes the proof.\hfill$\Box$


\section{Monotonicity with respect to $y$}\label{monotonocity section}

This section is devoted to the proof of Theorem \ref{phi is y increasing}. To furnish this goal, we need to establish a generalized comparison principle in unbounded domains of the form $C_{\alpha,\beta,l}^+.$ Then, together with further estimates on the behavior of any solution $\phi$ of the problem (\ref{c equation})-(\ref{cc}) in the lower cone $C_{\alpha,\beta,l}^-$ and with some ``sliding techniques'' which are similar to those done by Berestycki and Nirenberg \cite{c3bn91}, we prove that the solution $\phi$ is increasing in~$y$.

Let us first state the following proposition which is an important step to prove the main result in this section.

\begin{proposition}\label{liminf d_y phi>0}
Let $\alpha$ and $\beta$ belong to $(0,\pi)$. If $(c,\phi)$ is a solution of $(\ref{c equation})$-$(\ref{cc})$, then
$$\displaystyle{\Lambda:=\liminf_{l\rightarrow-\infty}\Big(\inf_{(x,y)\in C_{\alpha,\beta,l}^-}\frac{\partial_y\phi(x,y)}{\phi(x,y)}\Big)>0.}$$
\end{proposition}

\noindent\textbf{Proof.} Similar to the discussion in \cite{c3bh02}, we get from standard Schauder interior estimates and Harnack inequalities that there exists a constant $K$ such that
\begin{equation}\label{schauder harnack}
\forall (x,y)\in\mathbb R^2,~~\left|\partial_y\phi(x,y)\right|\leq K\phi(x,y)\hbox{ and }\left|\partial_x\phi(x,y)\right|\leq K\phi(x,y).
\end{equation}
Consequently, the function $\displaystyle{{\partial_y\phi}/{\phi}}$ is globally bounded in $\mathbb R^2.$ Denote by
$$\displaystyle{\Lambda:=\liminf_{l\rightarrow-\infty}\inf_{(x,y)\in C_{\alpha,\beta,l}^-}\frac{\partial_y\phi(x,y)}{\phi(x,y)}}$$
and let $\{l_n\}_{n\in\mathbb{N}}$ and $\{(x_n,y_n)\}_{n\in\mathbb{N}}$ be two sequences such that $(x_n,y_n)\in C_{\alpha,\beta,l_n}^-$ for all $n\in \mathbb{N}$, $l_n\rightarrow-\infty$ as $n\rightarrow+\infty$, and
$$\displaystyle{\frac{\partial_y\phi(x_n,y_n)}{\phi(x_n,y_n)}}\rightarrow \Lambda\ \hbox{as}\ n\rightarrow+\infty.$$
Next, we will proceed in several steps to prove that $\Lambda >0$.\hfill\break

\noindent{\it Step 1: From $(\ref{c equation})$ to a linear elliptic equation}. For each $n\in\mathbb{N}$, let
$$\displaystyle{\phi^{n}(x,y)=\frac{\phi(x+x_n,y+y_n)}{\phi(x_n,y_n)}\quad \hbox{for all}\quad (x,y)\in\mathbb R^2}.$$
Owing to the equation (\ref{c equation}) satisfied by $\phi$, we know that each function $\phi^n(x,y)$ satisfies the following equation
$$\Delta\phi^{n}(x,y)+(q(x+x_n)-c)\partial_y\phi^n(x,y)+\frac{f(\phi(x+x_n,y+y_n))}{\phi(x+x_n,y+y_n)}\phi^{n}(x,y)=0$$
for all $(x,y)\in\mathbb R^2$. Moreover, for any given $(x,y)\in\mathbb R^2$, it follows from (\ref{cc}) that the sequence $\phi(x+x_n,y+y_n)\rightarrow0$ as $n\rightarrow+\infty$ (since $(x_n,y_n)\in C_{\alpha,\beta,l_n}^-$ for each $n\in\mathbb{N}$ and $l_n\rightarrow-\infty$ as $n\rightarrow+\infty$). Noticing that $f(0)=0$, then we have
$$\frac{f(\phi(x+x_n,y+y_n))}{\phi(x+x_n,y+y_n)}\rightarrow f'(0)$$
as $n\rightarrow+\infty$. Since the function $q$ is $L-$periodic, we can construct a sequence $\{\tilde{x}_n\}_{n\in\mathbb N}$ such that $\tilde{x}_n\in [0,L]$ for all $n\in\mathbb{N}$ and
$$\forall n\in \mathbb{N},~\forall x\in\mathbb R,~q_n(x):=q(x+x_n)=q(x+\tilde{x}_n).$$
Consequently, there exists a point $x_\infty\in [0,L]$ such that $\tilde{x}_n\rightarrow x_\infty$ as $n\rightarrow+\infty$ (up to extraction of some subsequence), and the functions $q_n(x)$ converge uniformly to $q(x+x_\infty)$. Observe also that the functions $\phi^n$ are locally bounded in $\mathbb R^2$, from the estimates (\ref{schauder harnack}). From the standard elliptic estimates, the functions $\phi^n$ converge in all $W^{2,p}_{loc}(\mathbb R^2)$ weak (for $1<p<\infty$), up to extraction of another subsequence, to a nonnegative function $\phi^{\infty}$ which satisfies the following linear elliptic equation
\begin{equation}\label{eq by phi^infty}
\Delta\phi^{\infty}+(q(x+x_\infty)-c)\partial_y\phi^\infty+f'(0)\phi^\infty=0\ \hbox{in}\ \mathbb R^2.
\end{equation}
Furthermore, by the definition of $\phi^n$, we have $\phi^\infty(0,0)=1$. Then, the strong maximum principle yields that the function $\phi^{\infty}$ is positive everywhere in $\mathbb R^2.$\hfill\break

\noindent{\it Step 2: The form of $\phi^\infty$}. For any given $(x,y)\in\mathbb R^2,$ we have
\begin{equation}\label{c3m}
\partial_y{\phi^n}(x,y)=\displaystyle{\frac{\partial_y\phi(x+x_n,y+y_n)}{\phi(x_n,y_n)}}=\frac{\partial_y\phi(x+x_n,y+y_n)}{\phi(x+x_n,y+y_n)}\times\phi^n(x,y)
\end{equation}
for all $n\in\mathbb N.$ Referring to the definition of $\Lambda$, one can then conclude that for any given $(x,y)\in\mathbb R^2$,
$$\liminf_{n\rightarrow+\infty}\frac{\partial_y\phi(x+x_n,y+y_n)}{\phi(x+x_n,y+y_n)}\geq\Lambda.$$
Passing to the limit as $n\rightarrow+\infty$ in (\ref{c3m}) leads to
\begin{equation}\label{c3phiinfty1}
\partial_y\phi^{\infty}(x,y)\geq\Lambda\phi^{\infty}(x,y),\hbox{ for all }(x,y)\in\mathbb R^2.
\end{equation}
Furthermore,
\begin{equation}\label{c3phiinfty2}
\partial_y\phi^{\infty}(0,0)=\lim_{n\rightarrow+\infty}\partial_y\phi^n(0,0)=\lim_{n\rightarrow+\infty}\frac{\partial_y\phi(x_n,y_n)}{\phi(x_n,y_n)}=\Lambda=\Lambda\,\phi^{\infty}(0,0).
\end{equation}
Set
$$z^{\infty}(x,y)=\frac{\partial_y\phi^\infty(x,y)}{\phi^\infty(x,y)}\ \hbox{for all}\ (x,y)\in\mathbb R^2.$$
The function $z^{\infty}(x,y)$ is then a classical solution of the equation
\begin{equation}\label{eq by z^infty}
\Delta z^{\infty}+w\cdot\nabla z^{\infty}=0~\hbox{ in }~\mathbb R^2,
\end{equation}
where
$$\displaystyle{w=w(x,y)=\left(2\frac{\partial_x\phi^{\infty}}{\phi^\infty},2\frac{\partial_y\phi^{\infty}}{\phi^\infty}+q(x+x_\infty)-c\right)}$$
is a globally bounded vector field defined in $\mathbb R^2$ (see (\ref{schauder harnack})). It follows from (\ref{c3phiinfty1}) and (\ref{c3phiinfty2}) that
$$z^\infty(0,0)=\Lambda \quad \hbox{and}\quad z^\infty(x,y)\geq\Lambda \ \hbox{for all}\ (x,y)\in\mathbb R^2.$$
Obviously, the constant function $\Lambda$ also solves (\ref{eq by z^infty}). Then, it follows from the strong maximum principle that $z^\infty(x,y)=\Lambda$ for all $(x,y)\in\mathbb R^2$, and thus,
$$\forall (x,y)\in\mathbb R^2,~\phi^{\infty}(x,y)=e^{\Lambda y}\psi(x)>0$$
for some positive function $\psi(x)$ defined in $\mathbb R.$ Owing to (\ref{eq by phi^infty}), the function $\psi(x)$ is then a classical solution of the following ordinary differential equation
\begin{equation}\label{eq by psi(x)}
\psi''(x)+\left(\Lambda^2+\Lambda \,q(x+x_\infty)-c\Lambda+f'(0)\right)\psi(x)=0 \ \hbox{for all}\ x\in\mathbb R.
\end{equation}

\noindent{\it Step 3: From $(\ref{eq by psi(x)})$ to an eigenvalue problem}. Let
$$\displaystyle{\mu=\inf_{x\in\mathbb R}\frac{\psi(x+L)}{\psi(x)}},$$
where $L$ is the period of $q$ (see (\ref{cq})). From (\ref{schauder harnack}), the function $\psi$ satisfies $|\psi'(x)|\le|K\psi(x)|$ for all $x\in\mathbb R$ and $\mu$ is then a real number. Let $\{x'_n\}_{n\in\mathbb N}$ be a sequence in $\mathbb R$ such that
$$\frac{\psi(x'_n+L)}{\psi(x'_n)}\rightarrow\mu\ \hbox{ as }n\rightarrow+\infty.$$
Define a sequence of functions $\{\psi^n(x)\}_{n\in\mathbb{N}}$ by
$$\psi^n(x)=\frac{\psi(x+x'_n)}{\psi(x'_n)}\ \hbox{for all}\ x\in\mathbb R.$$
Then, for each $n\in\mathbb{N}$, the function $\psi^n(x)$ satisfies
$$(\psi^{n})''(x)+\left(\Lambda^2+\Lambda \,q(x+x'_n+x_\infty)-c\Lambda+f'(0)\right)\psi^n(x)=0$$
for all $x\in\mathbb R$.

Similar to the discussion in Step 1 and also due to the $L-$periodicity of $q$, it easily follows that, up to extraction of a subsequence, $\psi^n\rightarrow\psi^\infty$ in $C^2_{loc}(\mathbb R^2)$ and
$$q(\cdot+x'_n+x_\infty)\rightarrow q(\cdot+x'_\infty)\ \hbox{ as }n\rightarrow+\infty\hbox{, uniformly on each compact of }\mathbb R,$$
for some $x'_{\infty}\in\mathbb R$. Furthermore, the function $\psi^\infty$ is a nonnegative classical solution of the following equation
\begin{equation}\label{eq by psi^infty}
(\psi^{\infty})''+\left(\Lambda^2+\Lambda \,q(x+x'_\infty)-c\Lambda+f'(0)\right)\psi^{\infty}=0~\hbox{ in }~\mathbb R.
\end{equation}
Since $\psi^n(0)=1$ for all $n\in\mathbb{N}$, we have $\psi^\infty(0)=1$. Then, the strong maximum principle yields that $\psi^\infty(x)>0$ for all $x\in\mathbb R$.

Now, we consider a new function
$$\displaystyle{h(x):=\frac{\psi^\infty(x+L)}{\psi^\infty(x)}},$$
which is defined in $\mathbb R$. By the definition of $\mu$ and $\psi^n$, we have
$$\frac{\psi^n(x+L)}{\psi^n(x)}=\frac{\psi(x+x'_n+L)}{\psi(x+x'_n)}\geq\mu,\hbox{ for all } n\in\mathbb N \hbox{ and } x\in\mathbb R.$$
Passing to the limit as $n\rightarrow+\infty,$ one gets $h(x)\geq\mu$ for all $x\in\mathbb R.$ Moreover,
$$\psi^\infty(L)=\lim_{n\rightarrow+\infty}\psi^{n}(L)=\lim_{n\rightarrow+\infty}\frac{\psi(x'_n+L)}{\psi(x'_n)}=\mu.$$
Denote by
$$v(x)=\psi^{\infty}(x+L)-\mu\,\psi^\infty(x)\ \hbox{for all}\ x\in\mathbb R.$$
Then, the function $v$ is nonnegative and satisfies the linear elliptic equation (\ref{eq by psi^infty}) with the property $v(0)=0$. Thus, the strong maximum principle yieldsthat $v\equiv 0$ in $\mathbb R$, and consequently, $h(x)=\mu>0$ in $\mathbb R$ (since $\psi^\infty(x)>0$ for all $x\in\mathbb R$).

Define $\theta=L^{-1}\ln\mu$. If we write $\psi^{\infty}(x)=e^{\theta x}\varphi(x)$ for all $x\in\mathbb R$, then it follows from $\psi^{\infty}(x+L)=\mu\,\psi^{\infty}(x)$ that
$$\forall x\in\mathbb R,~~\varphi(x+L)=\varphi(x).$$
After replacing  $\psi^{\infty}$ by $e^{\theta x}\varphi$ in (\ref{eq by psi^infty}), we conclude that the function $\varphi$ is a classical solution of the following problem
\begin{equation}\label{eq by varphi}\left\{\begin{array}{l}
\varphi''+2\theta\varphi'+\theta^2\varphi+\left(\Lambda^2-c\Lambda+q(x+x'_\infty)\Lambda+f'(0)\right)\varphi=0\hbox{ in }\mathbb R, \vspace {4 pt}\\
\varphi \hbox{ is $L$-periodic},\vspace {4 pt}\\
\forall x\in\mathbb R,~\varphi(x) >0.\end{array}\right.
\end{equation}
For each $\lambda\in\mathbb R$, we define an elliptic operator as follows
$$L_{\theta,\lambda}:=\frac{d^2}{dx^2}+2\theta\frac{d}{dx}+\left[\theta^2+\lambda^2-c\lambda+q(x+x'_\infty)\lambda+f'(0)\right]$$
acting on the set
$$E:=\{g(x)\in C^2(\mathbb{R}); g(x+L)=g(x)\ \hbox{for all}\ x\in\mathbb R \}.$$
We denote by $k_\theta(\lambda)$ and $\varphi^{\theta,\lambda}$ the principal eigenvalue and the corresponding principal eigenfunction of this operator. In addition to the existence, we also have the uniqueness (up to a multiplication by any nonzero constant) of the principal eigenfunction $\varphi^{\theta,\lambda}$ which keeps sign over $\mathbb R$ and solves the following problem
\begin{equation}\label{eigenvalue pb}\left\{\begin{array}{l}
L_{\theta,\lambda}\varphi^{\theta,\lambda}=k_\theta(\lambda)\varphi^{\theta,\lambda}\hbox{ in }\mathbb R\vspace{4 pt}\\
\varphi^{\theta,\lambda} \hbox{ is $L-$periodic}.\end{array}\right.
\end{equation}

From (\ref{eq by varphi}) and the above discussions, we conclude that, for $\lambda=\Lambda$, $k_\theta(\Lambda)=0$ is the principal eigenvalue and the function $\varphi$ is the corresponding eigenfunction. In other words,~$\Lambda$ is a solution of the equation $k_\theta(\lambda)=0$.

\begin{figure}
\begin{center}
\includegraphics[width=7 cm]{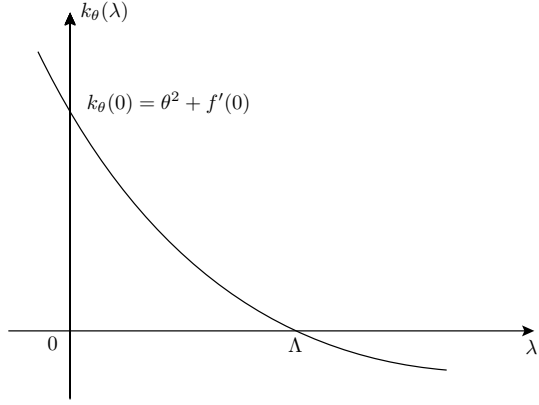}\\
\caption{The function $\lambda\mapsto k_\theta(\lambda).$}\label{eigenvalue}
\end{center}
\end{figure}

Now, we consider the function $\mathbb R\ni\lambda\mapsto k_\theta(\lambda)$. It follows from Proposition 5.7 in \cite{c3bh02} that $\lambda\mapsto k_\theta(\lambda)$ is of convex. Moreover, for $\lambda=0,$ the principal eigenfunction $\varphi^{\theta,0}$ is a constant function, say $\varphi^{\theta,0}\equiv 1$ (due to the uniqueness up to multiplication by a constant), and the principal eigenvalue is
$$k_{\theta}(0)=\theta^2+f'(0)>0.$$
Thus, in order to obtain that $\Lambda>0,$ it suffices to prove that $\displaystyle{\frac{d\,k_\theta}{d\lambda}(0)<0}$ (see figure \ref{eigenvalue}).

Since $\varphi^{\theta,\lambda}$ is $L-$periodic for each $\lambda\in\mathbb R$, we then integrate the equation (\ref{eigenvalue pb}) with respect to $x$ over $[0,L]$ to obtain
\begin{equation}\label{int by parts}\begin{array}{ll}
\displaystyle{k_\theta(\lambda)\int_{0}^L\varphi^{\theta,\lambda}(x)dx}=&\displaystyle{\underbrace{\left(\theta^2+f'(0)\right)}_{k_\theta(0)}\int_{0}^L\varphi^{\theta,\lambda}(x)dx+\lambda^2\int_{0}^L\varphi^{\theta,\lambda}(x)dx}\\
&\displaystyle{ -c\lambda\int_{0}^L\varphi^{\theta,\lambda}(x)dx+\lambda\int_{0}^Lq(x+x'_\infty)\varphi^{\theta,\lambda}(x)dx}\end{array}
\end{equation}
for all $\lambda\in\mathbb R$. Owing to standard elliptic estimates, the family $\{\varphi^{\theta,\lambda}\}_{\lambda\in\mathbb R}$, when normalized by $\max_{\mathbb R}\varphi^{\theta,\lambda}=1$, converges in $C^2_{loc}(\mathbb R)$ to the constant function $\varphi^{\theta,0}\equiv1$ as $\lambda$ converges to~$0$. Passing to the limit as $\lambda\rightarrow0$ in (\ref{int by parts}), one consequently gets
$$\displaystyle{\lim_{\lambda\rightarrow0}\frac{k_\theta(\lambda)-k_\theta(0)}{\lambda}=-c+\frac{1}{L}\int_{0}^Lq(x+x'_\infty)\,dx}.$$
However, by the assumptions (\ref{cq}) on $q$, we know that
$$\displaystyle{\int_{0}^L q(x+x'_\infty)\,dx=\int_{0}^L q(x)\,dx=0}.$$
Therefore,
$$\displaystyle{\frac{d k_\theta}{d\lambda}(0)=-c}.$$
But, from part ii) of Theorem \ref{c3theorem1}, the speed $c$ satisfies
$$c\ge c^*=\max \left(\frac{c^*_{A,q\sin\alpha,f}}{\sin\alpha},\frac{c^*_{B,q\sin\beta,f}}{\sin\beta}\right)>0.$$
Thus, $\frac{d k_\theta}{d\lambda}(0)<0$ and that completes the proof of Proposition \ref{liminf d_y phi>0}.\hfill$\Box$\break

In the following, we are going to establish a generalized comparison principle which will be an important tool in the proof of Theorem \ref{phi is y increasing}. Before stating this result, let us first introduce some notations and assumptions that we need in our setting. For each $l\in\mathbb R,$ $\alpha,\beta\in(0,\pi)$, we consider $A(x,y)=(A_{ij}(x,y))_{1\leq i,j\leq N}$ as a symmetric $C^{2,\delta}\big(\overline{C_{\alpha,\beta,l}^+}\big)$ matrix field satisfying
\begin{eqnarray}\label{cA}\left\{\begin{array}{l}
\exists 0<\alpha_1\leq\alpha_2, \; \forall(x,y)\in\overline{C_{\alpha,\beta,l}^+},\;\forall \xi\in\mathbb{R}^2, \vspace{4 pt}\\
\displaystyle{\alpha_1|\xi|^2 \leq\sum_{1\leq i,j\leq 2}A_{ij}(x,y)\xi_i\xi_j\leq\alpha_2|\xi|^2}.\end{array}\right.
\end{eqnarray}
Moreover,
$$\begin{array}{ll}
\partial C_{\alpha,\beta,l}^+:=&\Big\{(x,y)\in\mathbb R^2,~y=-x\cot\beta+l\hbox{ when }x\geq0,\vspace{3pt}\\
&\quad\quad\quad\quad\hbox{and }y=x\cot\alpha+l\hbox{ when }x\leq0\Big\}\end{array}$$
denotes the boundary of the subset $C_{\alpha,\beta,l}^+$ which was introduced in Definition \ref{cones}, and
$$\hbox{dist}\left((x,y);\partial C_{\alpha,\beta,l}^+\right)$$
stands for the Euclidean distance from $(x,y)\in\mathbb R^2$ to the boundary $\partial C_{\alpha,\beta,l}^+.$

The generalized comparison principle is now stated in the following lemma.

\begin{lemma}\label{comp principle} 
Let $\alpha$ and $\beta$ be fixed in $(0,\pi)$ and $l\in\mathbb{R}$. Let $g(x,y,u)$ be a globally bounded and a globally Lipschitz-continuous function defined in $\overline{C_{\alpha,\beta,l}^+}\times\mathbb{R}.$ Assume that $g$ is non-increasing with respect to $u$ in $\mathbb R^2\times[1-\rho,+\infty)$ for some $\rho>0.$ Let $\tilde{q}=\left(q_1(x,y),q_2(x,y)\right)$ be a globally bounded $C^{0,\delta}\big(\overline{C_{\alpha,\beta,l}^+}\big)$ vector field $($with $\delta>0)$ and let $A(x,y)=(A_{ij}(x,y))_{1\leq i,j\leq 2}$ be a symmetric $C^{2,\delta}\big(\overline{C_{\alpha,\beta,l}^+}\big)$ matrix field satisfying~$(\ref{cA})$.\par
Assume that $\phi^{1}(x,y)$ and $\phi^{2}(x,y)$ are two bounded uniformly continuous functions defined in $\overline{C_{\alpha,\beta,l}^+}$ of class $C^{2,\mu}\big(\overline{C_{\alpha,\beta,l}^+}\big)$ (for some $\mu>0$). Furthermore, we assume that
$$\left\{\begin{array}{ccc}
L\,\phi^{1}+g(x,y,\phi^{1})&\geq&0~\hbox{ in }~C_{\alpha,\beta,l}^+,\vspace{5 pt}\\
L\,\phi^{2}+g(x,y,\phi^{2})&\leq&0~\hbox{ in }~C_{\alpha,\beta,l}^+,\vspace{5 pt}\\
\phi^{1}(x,y)\leq \phi^{2}(x,y)&&~~\hbox{on }~\partial C_{\alpha,\beta,l}^+ ,\end{array}\right.$$
and that
\begin{equation}\label{as distance goes to infty}
\displaystyle{\limsup_{\displaystyle{(x,y)\in C_{\alpha,\beta,l}^+,\,{\rm dist}\left((x,y);\partial C_{\alpha,\beta,l}^+\right)\rightarrow+\infty}}\,[\phi^{1}(x,y)-\phi^{2}(x,y)]}\leq0,
\end{equation}
where $L$ is the elliptic operator defined by
$$L\phi:=\nabla_{x,y}\cdot(A\nabla_{x,y}\phi)+\tilde{q}(x,y)\cdot\nabla_{x,y}\phi.
$$\par
If $\phi^{2}\geq 1-\rho$ in $\overline{C_{\alpha,\beta,l}^+},$ then
$$\phi^{1}\leq\phi^{2}~~\hbox{in}~~\overline{C_{\alpha,\beta,l}^+}.$$
\end{lemma}

\begin{remark}{\rm Note here that $\phi^{1},\,\phi^{2},\,\tilde{q},\,A$ and $g$ are not assumed to be $L-$periodic with respect to $x.$}
\end{remark}

\noindent\textbf{Proof.} Since the functions $\phi^{1}$ and $\phi^{2}$ are globally bounded, one can then find $\varepsilon>0$ large enough such that $\phi^{1}-\varepsilon\leq\phi^{2}$ in $\overline{C_{\alpha,\beta,l}^+}$. Let us set
$$\displaystyle{\varepsilon^{*}=\inf\left\{\varepsilon>0,\;\phi^{1}-\varepsilon\leq\phi^{2}\;\hbox{in}\;\displaystyle{\overline{C_{\alpha,\beta,l}^+}}\right\}}\geq0.$$
By continuity, we then get $\phi^{1}-\varepsilon^{*}\leq\phi^{2}$ in $\displaystyle{\overline{C_{\alpha,\beta,l}^+}}.$ Thus, to complete the proof of Lemma~\ref{comp principle}, it suffices to prove that $\varepsilon^{*}=0.$

Assume $\varepsilon^{*}>0$. Then, there exist a sequence $\{\varepsilon_{n}\}_{n\in\mathbb N}$ converging to $\varepsilon^{*}$, with $0<\varepsilon_{n}<\varepsilon^*$ for all $n$, and a sequence of points $(x_n,y_n)\in\displaystyle{\overline{C_{\alpha,\beta,l}^+}}$ such that
$$\phi^{1}(x_n,y_n)-\varepsilon_n\geq\phi^{2}(x_n,y_n)\ \hbox{for all}\ n\in\mathbb{N}.$$

Because of (\ref{as distance goes to infty}) and since $\varepsilon^*>0$, the sequence $\displaystyle{\left\{{\rm dist}\left(\displaystyle{(x_n,y_n);\partial C_{\alpha,\beta,l}^+}\right)\right\}_{n\in\mathbb{N}}}$ is bounded. Furthermore, the facts that $\phi^{1}\leq \phi^{2}~\hbox{ on }~\partial C_{\alpha,\beta,l}^+ $ and $\phi^{1},$  $\phi^{2}$ are uniformly continuous yield that
$$R:=\displaystyle{\liminf_{n\rightarrow+\infty}\;{\rm dist}\left(\displaystyle{(x_n,y_n);\partial C_{\alpha,\beta,l}^+}\right)>0}.$$
For each $n\in\mathbb{N},$ let $(x'_n,y'_n)$ be a point on $\partial C_{\alpha,\beta,l}^+$ such that
$${\rm dist}\left(\displaystyle{(x_n,y_n);\partial C_{\alpha,\beta,l}^+}\right)=\left|(x'_n,y'_n)-(x_n,y_n)\right|.$$
Up to extraction of some subsequence, we can then conclude that there exists $(\bar x,\bar y)\in\mathbb R^2$ with $\left|(\overline{x},\overline{y})\right|=R$ such that
$$(x'_n,y'_n)-(x_n,y_n)\rightarrow(\overline{x},\overline{y})\hbox{ as } n\rightarrow+\infty.$$
Call $B_R:=\{(x,y)\in\mathbb R^2,~|(x,y)|<R\}$. It follows from the definition of $R$ that for any point $(x,y)\in B_R$ and for any $n\in\mathbb{N}$ large enough, we have $(x,y)+(x_n,y_n)\in C_{\alpha,\beta,l}^+$.

For each $(x,y)\in B_R$, call
$$\phi^{1}_n(x,y)=\phi^{1}(x+x_n,y+y_n)\hbox{ and }\phi^{2}_n(x,y)=\phi^{2}(x+x_n,y+y_n)$$
for $n$ large enough.

From the regularity assumptions on $\phi^{1}$ and $\phi^{2}$ and up to extraction of some subsequence, the functions $\phi^{i}_{n}$ converge in $C^{2}_{loc}(B_R)$ to two functions $\phi^{i}_{\infty}$ which can be extended by continuity to $\partial B_R$ and are of class $C^{2,\mu}\left(\overline{B_R}\right)$, for $i=1,2.$ Similarly, since $\tilde q$ and $A$ are globally $C^{0,\delta}\big(\overline{C_{\alpha,\beta,l}^+}\big)$ (for some $\delta>0$), we can assume that the fields $\tilde{q}_{n}(x,y)=\tilde{q}(x+x_n,y+y_n)$ and $A_n(x,y)=A(x+x_n,y+y_n)$ converge as $n\rightarrow+\infty$ in $B_R$ to two fields $\tilde{q}_{\infty}$ and $A_{\infty}$ which are of class $C^{0,\delta}\big(\overline{B_R}\big)$. The matrix $A_\infty$ satisfies the same ellipticity condition as $A$ which is given in (\ref{cA}).

For each $(x,y)\in B_R$, the functions $\phi^{i}_{n},$ $i=1,2,$ satisfy
$$\displaystyle {L_n\,\phi^{1}_{n}-L_n\,\phi^{2}_{n}\geq-g(x+x_n,y+y_n,\phi_{n}^{1}(x,y))+g(x+x_n,y+y_n,\phi^{2}_{n}(x,y))}$$
for $n$ large enough, where
$$L_n\phi:=\displaystyle{\nabla_{x,y}\cdot(A_{n}\nabla_{x,y}\phi)}\displaystyle{\,+\,\tilde{q}_{n}\cdot\nabla_{x,y}\phi.}$$
Since $\phi^{2}\geq 1-\rho$ in $\overline{C_{\alpha,\beta,l}^+}$ and $g(x,y,u)$ is non-increasing with respect to $u$ in the set $\overline{C_{\alpha,\beta,l}^+}\times[1-\rho,+\infty),$ we get
\begin{equation}\label{L_n (phi_n)^1 -(phi_n)^{2}}
\begin{split}
L_n\,\phi^{1}_{n}-L_n\,\phi^{2}_{n}\geq&-g(x+x_n,y+y_n,\phi^{1}_{n}(x,y))\\
&+g(x+x_n,y+y_n,\phi^{2}_{n}(x,y)+\varepsilon^{*}).
\end{split}
\end{equation}

From the assumptions of Lemma \ref{comp principle}, we can also assume, up to extraction of some subsequence, that the functions
$$\displaystyle{R_{n}(x,y):=-g(x+x_n,y+y_n,\phi^{1}_{n}(x,y))+g(x+x_n,y+y_n,\phi^{2}_{n}(x,y)+\varepsilon^{*})}$$
converge to a function $R_{\infty}(x,y)$ locally uniformly in $B_R$. Since
$$|R_n(x,y)|\leq||g||_{Lip}|\phi^{1}_{n}(x,y)-\varepsilon^{*}-\phi^{2}_{n}(x,y)|$$
for all $n\in\mathbb N$, we get $|R_\infty(x,y)|\leq||g||_{Lip}|\phi^{1}_{\infty}(x,y)-\varepsilon^{*}-\phi^{2}_{\infty}(x,y)|$. In other words, there exists a globally bounded function $B(x,y)$ defined in $B_R$ such that
$$R_{\infty}(x,y)=B(x,y)\left[\phi^{1}_{\infty}(x,y)-\varepsilon^*-\phi^{2}_{\infty}(x,y)\right]\ \hbox{for all}\ (x,y)\in B_R.$$
By passing to the limit as $n\rightarrow+\infty$ in (\ref{L_n (phi_n)^1 -(phi_n)^{2}}), it follows that
$$L_\infty\phi^{1}_\infty-L_\infty\phi^{2}_\infty\geq B(x,y)(\phi^{1}_\infty-\varepsilon^{*}-\phi^{2}_\infty)\hbox{ in }B_R,$$
where $L_\infty\phi:=\displaystyle{\nabla_{x,y}\cdot(A_{\infty}\nabla_{x,y}\phi)}\displaystyle{\,+\,\tilde{q}_{\infty}\cdot\nabla_{x,y}\phi.}$ Let
$$z(x,y)=\phi^{1}_{\infty}-\varepsilon^{*}-\phi^{2}_{\infty}\hbox{ in }\overline{B_R}.$$
We then get
\begin{equation}\label{equation by z}
L_\infty z-B(x,y)z\geq0\hbox{ in }B_R.
\end{equation}
Noticing that $(x'_n,y'_n)\in \partial C_{\alpha,\beta,l}^+$, that $\phi^1\leq\phi^2$ over $\partial C_{\alpha,\beta,l}^+$, that $\phi^1$ and $\phi^2$ are uniformly continuous in $\overline{C^+_{\alpha,\beta,l}}$, and that $(x'_n,y'_n)-(x_n,y_n)\rightarrow (\overline{x},\overline{y})$, we have
\begin{equation}\label{on boundary of B_R}
\phi^1_\infty(\overline{x},\overline{y})\leq \phi^2_\infty(\overline{x},\overline{y}).
\end{equation}
On the other hand, for each $(x,y)\in B_R$, $\phi^1_n(x,y)-\varepsilon^*\leq\phi^2_n(x,y)$ for $n$ large enough, and $\phi^1_n(0,0)-\varepsilon_n\geq \phi^2_n(0,0)$. Passing to the limit as $n\rightarrow+\infty$ and over $\partial B_R$, then by continuity, we get
$$\phi^1_\infty(x,y)-\varepsilon^*\leq\phi^2_\infty(x,y)\hbox{ in }\overline{B_R},$$
and
$$\phi^1_\infty(0,0)-\varepsilon^*=\phi^2_\infty(0,0).$$
Consequently, the function $z=z(x,y)$ is a nonpositive continuous function in $\overline{B_R}$, satisfying~(\ref{equation by z}) in $B_R$ and such that $z(0,0)=0$. Then, the strong maximum principle yields that $z\equiv 0$ in $\overline{B_R}$ with $\varepsilon^*>0.$ Namely, $\phi^1_\infty(x,y)-\varepsilon^*=\phi^2_\infty(x,y)$ for all $(x,y)\in\overline{B_R}$. We get a contradiction with (\ref{on boundary of B_R}) by choosing $(x,y)=(\overline{x},\overline{y})$ ($\in\partial B_R$).\hfill$\Box$\break

The following lemma is devoted to proving the positivity of the infimum of a conical front solving (\ref{c equation}-\ref{cc}) over any set having the form of an ``upper cone''. This lemma will be also used in the proof of Theorem \ref{phi is y increasing}.

\begin{lemma}\label{inf is positive over upp cones}
For any fixed $\alpha$ and $\beta$ in $(0,\pi)$, let $(c,\phi)$ be a solution of $(\ref{c equation})$-$(\ref{cc})$. Then,
\begin{equation}\label{inf phi >0 on upp cones}
\forall\, l\in\mathbb R,~~\inf_{(x,y)\in C_{\alpha,\beta,l}^+}\phi(x,y)>0.
\end{equation}
\end{lemma}

\noindent\textbf{Proof.} Since the function $\phi$ is nonnegative in $\mathbb R^2$, then $\inf_{\mathbb R^2}\phi\geq0$. In order to prove (\ref{inf phi >0 on upp cones}), we assume to the contrary that $\inf_{(x,y)\in C_{\alpha,\beta,l_0}^+}\phi(x,y)=0$ for some fixed $l_0\in\mathbb R$. Thus, there exists a sequence $\{(x_n,y_n)\}_{n\in\mathbb N}$ in $C_{\alpha,\beta,l_0}^+$ such that $\phi(x_n,y_n)\rightarrow0$ as $n\rightarrow+\infty.$ On the other hand, the limiting condition $\displaystyle{\lim_{l\rightarrow+\infty}\inf_{(x,y)\in C^{+}_{\alpha,\beta,l}}\phi(x,y)= 1}$ yields that there exists $M\in\mathbb R$ such that
\begin{equation}\label{geq 3/4}
\displaystyle{\forall (x,y)\in C_{\alpha,\beta,M}^+,~~ \phi(x,y)\geq\frac{3}{4}}.
\end{equation}
We recall that $\displaystyle{{\rm dist}\left((x_n,y_n);\partial C^{+}_{\alpha,\beta,l_0}\right)}$ is the Euclidean distance from $(x_n,y_n)\in\mathbb R^2$ to the
boundary $\partial C_{\alpha,\beta,l_0}^+$. Having (\ref{geq 3/4}) and the fact that $\phi(x_n,y_n)\rightarrow0$ as $n\rightarrow+\infty$, we know that the sequence $\displaystyle{\{\hbox{dist}((x_n,y_n);\partial C^{+}_{\alpha,\beta,l_0})\}_{n\in\mathbb N}}$ should be bounded and consequently,
\begin{equation}\label{translate to high cone}
\exists\,(\overline{x},\overline{y})\in\mathbb R^2\hbox{ such that }(\overline{x}+x_n,\overline{y}+y_n)\in C_{\alpha,\beta,M}^+
\end{equation}
for all $n\in\mathbb N$. Now, we define $\phi_n(x,y):=\phi(x+x_n,y+y_n)$ for all $(x,y)\in\mathbb R^2$ and $n\in\mathbb{N}$. From (\ref{c equation}), the function $\phi_n$ is a classical solution of the following equation
$$\Delta_{x,y}\phi_n+(q(x+x_n)-c)\partial_y\phi_n+f(\phi_n)=0\hbox{ in }\mathbb R^2,$$
for all $n\in\mathbb N.$

The function $q$ is a globally bounded $C^{0,\delta}(\mathbb R)$ function which is $L-$periodic. As a consequence, we can assume that the sequence of functions $q_n(x):=q(x+x_n)$ converges uniformly in $\mathbb R$, as $n\rightarrow+\infty$, to the function $q_\infty:=q(x+x_\infty)$ for some $x_\infty\in\mathbb R$. The regularity of the function $\phi$ yields that the sequence $\{\phi_n\}_{n\in\mathbb N}$ is bounded in $C^{2,\delta}(\mathbb R^2)$. Thus, up to extraction of some subsequence, $\phi_n\rightarrow\phi_\infty$ in $C^2_{loc}(\mathbb R^2)$ as $n\rightarrow+\infty$, where $\phi_\infty$ is a nonnegative ($0\leq\phi_n\leq1$ for all $n\in\mathbb N$) classical solution of the equation
$$\Delta_{x,y}\phi_\infty+\left(q(x+x_\infty)-c\right)\partial_y\phi_\infty+f(\phi_\infty)=0\hbox{ in }\mathbb R^2.$$
Moreover, $\phi_\infty(0,0)=\lim_{n\rightarrow+\infty}\phi(x_n,y_n)=0$.

Since $f\geq0$ in $[0,1]$, we then have
$$\left\{\begin{array}{l}
\Delta_{x,y}\phi_\infty+\left(q(x+x_\infty)-c\right)\partial_y\phi_\infty\leq0\hbox{ in }\mathbb R^2,\vspace{3 pt}\\
0\leq\phi_\infty\leq 1\hbox{ in }\mathbb R^2,\vspace{3 pt}\\
\phi_\infty(0,0)=0.\end{array}\right.$$
The strong maximum principle implies that $\phi_\infty\equiv 0$ in $\mathbb R^2$. However, we can conclude from~(\ref{geq 3/4}) and~(\ref{translate to high cone}) that
$$\forall n\in\mathbb N,~\phi(\overline{x}+x_n,\overline{y}+y_n)\geq \frac{3}{4}.$$
Passing to the limit as $n\rightarrow+\infty,$ one gets $\phi_\infty(\overline{x},\overline{y})\geq {3}/{4}$, which is a contradiction with $\phi_\infty\equiv0$ in $\mathbb R^2.$ Therefore, our assumption that $\inf_{(x,y)\in C_{\alpha,\beta,l_0}^+}\phi(x,y)=0$ is false and that completes the proof of Lemma \ref{inf is positive over upp cones}.\hfill $\Box$\break

Now, we are in the position to give the proof of the main result in this section.\hfill\break

\noindent\textbf{Proof of Theorem \ref{phi is y increasing}.} In this proof, we call
$$\forall \tau \in\mathbb R,~\phi^{\tau}(x,y):=\phi(x,y+\tau)\ \hbox{for all}\ (x,y)\in\mathbb R^{2}.$$
Assume that one has proved that $\phi^\tau\geq\phi$ in $\mathbb R^2$ for all $\tau\geq0.$ Since the coefficients $q$ and $f$ are independent of $y$, then for any $h>0$, the nonnegative function $z(x,y):=\phi^h(x,y)-\phi(x,y)$ is a classical solution (due to (\ref{c equation})) of the following linear elliptic equation
$$\Delta_{x,y}z+(q(x)-c)\partial_yz+b(x,y)z=0\hbox{ in }\mathbb R^2,$$
for some globally bounded function $b=b(x,y)$. It follows from the strong maximum principle that the function $z$ is either identically $0,$ or positive everywhere in $\mathbb R^2$. Due to the conical limiting conditions (\ref{cc}) satisfied by the function $\phi$, we can conclude that the function $z$ can not be identically $0$. In fact, if $z\equiv0$, then $\phi(x,y+h)=\phi(x,y)$ for all $(x,y)\in\mathbb R^2$ with $h>0.$ This yields that $\phi$ is $h-$periodic with respect to $y$, which is impossible from (\ref{cc}). Hence, the function $z$ is positive everywhere in $\mathbb R^2$, and consequently, the function $\phi$ is increasing in~$y$.

By virtue of the above discussion, we only need to prove that $\phi^\tau\geq\phi$ for all $\tau\geq0$. Proposition \ref{liminf d_y phi>0} yields that there exists $l_0\in\mathbb R$ such that $\partial_y\phi(x,y)>0$ for all $(x,y)\in C_{\alpha,\beta,l_0}^-$. On the other hand, Lemma \ref{inf is positive over upp cones} yields that $\displaystyle{\inf_{(x,y)\in C^{+}_{\alpha,\beta,l_0}}\phi(x,y)}>0.$ Since
$$\displaystyle{\lim_{l\rightarrow-\infty}\sup_{(x,y)\in C^{-}_{\alpha,\beta,l}}\phi(x,y)= 0,}$$
there exists then $B>0$ such that $-B\leq l_0$ and
$$\forall (x,y)\in C_{\alpha,\beta,-B}^-\;,~~\phi(x,y)\leq \inf_{(x',y')\in C^{+}_{\alpha,\beta,l_0}}\phi(x',y'),$$
and consequently, we have
\begin{equation}\label{on lower cone}
\forall\,\tau\geq0,~\forall (x,y)\in C_{\alpha,\beta,-B}^-\;,~~ \phi(x,y)\leq \phi(x,y+\tau).
\end{equation}
The above inequality is indeed satisfied in both cases $y+\tau\le l_0$ and $y+\tau\ge l_0$. The assumption that $f'(1)<0$ in (\ref{cf}) and the continuity of $f'$ over $[0,1]$ lead to the existence of $0<\eta<1$ such that $f$ is non-increasing in $[1-\eta,1]$. Furthermore, even it means increasing $B,$ one can assume, due to (\ref{cc}), that $\phi(x,y)\geq 1-\eta$ for all $(x,y)\in \overline{C_{\alpha,\beta,B}^+}$ and $\phi(x,y)\leq \theta$ for all $(x,y)\in\overline{C_{\alpha,\beta,-B}^-}$, where $\theta$ is choosen so that $0<\theta<1-\eta$. We apply Lemma \ref{comp principle} to the functions $\phi^1:=\phi$ and $\phi^2:=\phi^\tau$ with $\tau\ge 2B$, by taking $\rho=\eta$, $A=I$, $g=f$, $\tilde{q}(x)=\left(0,q(x)-c\right)$ in $\mathbb R$ and $l\,=\,-B,$ to obtain
$$\forall \tau\geq2B,~\forall(x,y)\in\overline{C_{\alpha,\beta,-B}^+}\;,\;\phi(x,y)\leq\phi^{\tau}(x,y).$$
Combining the above inequality with (\ref{on lower cone}), we have
$$\forall \tau\geq2B,~\forall(x,y)\in \mathbb R^2,\ \phi(x,y)\leq\phi^{\tau}(x,y).$$

Let us now decrease $\tau$ and set
$$\tau^{*}=\inf\left\{\tau>0,\phi(x,y)\leq\phi(x,y+\tau')\ \hbox{for all}\ \tau'\geq\tau\ \hbox{and for all}\ (x,y)\in\mathbb R^2\,\right\}.$$
First, we note that $\tau^{*}\leq2B$, and by continuity, we have $\phi\leq\phi^{\tau^*}~\hbox{in}~\mathbb{R}^2.$ Call
$$S:=C_{\alpha,\beta,-B}^+\setminus C_{\alpha,\beta,B}^{-}$$
the slice located between the ``lower cone'' $C_{\alpha,\beta,-B}^{-}$ and the ``upper cone'' $C_{\alpha,\beta,B}^+$. Then, for the value of $\displaystyle{\sup_{(x,y)\in \overline{S}}\left(\phi(x,y)-\phi^{\tau^*}(x,y)\right)}$, the following two cases may occur.\hfill\break

\noindent{\textit{Case 1:}} suppose that
$$\displaystyle{\sup_{(x,y)\in \overline{S}}\left(\phi(x,y)-\phi^{\tau^*}(x,y)\right)<0}.$$
Since the function $\phi$ is (at least) uniformly continuous, there exists $\varepsilon>0$ such that $0<\varepsilon<\tau^*$ and the above inequality holds for all $\tau\in[\tau^*-\varepsilon,\tau^*]$. Then, for any $\tau$ in the interval~$[\tau^*-\varepsilon,\tau^*]$, due to (\ref{on lower cone}) and the definition of $S$, we get that
$$\phi(x,y)\leq \phi^\tau(x,y)\hbox{ over }\overline{C_{\alpha,\beta,B}^{-}}.$$
Hence, $\phi\leq \phi^\tau$ over $\partial C_{\alpha,\beta,B}^{+}$. On the other hand, since $\tau\geq \tau^*-\varepsilon>0$ and $\phi\geq 1-\eta$ over~$\overline{C_{\alpha,\beta,B}^{+}}$, we have $\phi^\tau\geq1-\eta$ over $\overline{C_{\alpha,\beta,B}^{+}}$. Lemma~\ref{comp principle},  applied to $\phi$ and $\phi^{\tau}$ in $C_{\alpha,\beta,B}^{+},$ yields that
$$\phi(x,y)\leq\phi^{\tau}(x,y)\ \hbox{for all}\ (x,y)\in\overline{ C_{\alpha,\beta,B}^{+}}.$$
As a consequence, we obtain $\phi\leq\phi^{\tau}$ in $\mathbb{R}^2,$ and that contradicts the minimality of $\tau^{*}.$ Therefore, case 1 is ruled out.\hfill\break

\noindent{\textit{Case 2:}} suppose that
$$\displaystyle{\sup_{(x,y)\in \overline{S}}\left(\phi(x,y)-\phi^{\tau^*}(x,y)\right)=0.}$$
Then, there exists a sequence of points $\{(x_n,y_n)\}_{n\in\mathbb N}$ in $\overline{S}$ such that
\begin{equation}\label{at x_n,y_n}
\phi(x_n,y_n)-\phi^{\tau^*}(x_n,y_n)\rightarrow 0~\hbox{as}~n\rightarrow+\infty.
\end{equation}

For each $n\in\mathbb{N},$ call $\phi_n(x,y)=\phi(x+x_n,y+y_n)\hbox{ and }\phi^{\tau^*}_n(x,y)=\phi^{\tau^*}(x+x_n,y+y_n),$ for all $(x,y)\in \mathbb R^2.$ From the regularity assumptions for $\phi$ and up to extraction of some subsequence, the functions $\phi_{n}$ and $\phi^{\tau^*}_n$ converge in $C^{2}_{loc}(\mathbb R^2)$ to two functions $\phi_{\infty}$ and $\phi^{\tau^*}_\infty$ in~$C^{2,\delta}(\mathbb R^2).$ On the other hand, since $q$ is globally $C^{0,\delta}\left(\mathbb R\right)$ and $L-$periodic, we can assume that the functions $q_{n}(x)=q(x+x_n)$ converge locally in $\mathbb R$ to a globally $C^{0,\delta}\left(\mathbb R\right)$ function~$q_{\infty}$ as $n\rightarrow+\infty.$

For any $(x,y)\in\mathbb R^2$, set $z(x,y)=\phi_\infty(x,y)-\phi^{\tau^{*}}_\infty(x,y)$. The function $z$ is nonpositive because $\phi\leq\phi^{\tau^*}$ in $\mathbb R^{2}$. Moreover, by passing to the limit as $n\rightarrow+\infty$ in (\ref{at x_n,y_n}), we obtain $z(0,0)=0$. Furthermore, since the function $q$ does not depend on $y$, we know that the function $z$ solves the following linear elliptic equation
$$\Delta_{x,y}z+(q_{\infty}(x)-c)\partial_{y}z+b(x,y) z=0\hbox{ in }\mathbb R^2$$
for some globally bounded function $b(x,y)$ (since $f$ is Lipschitz continuous). Then, the strong elliptic maximum principle implies that either $z>0$ in $\mathbb R^2$ or $z=0$ everywhere in~$\mathbb R^2$. In fact, the latter case is impossible because it contradicts with the conical conditions at infinity (\ref{cc}): indeed, since $(x_n,y_n)\in \bar S$ for all $n\in\mathbb N$, it follows from (\ref{cc}) that $\lim_{y\rightarrow+\infty}\phi^{\infty}(0,y)=1$ and $\lim_{y\rightarrow-\infty}\phi^{\infty}(0,y)=0$, whence the function $\phi^{\infty}$ cannot be $\tau^*$-periodic with respect to $y$, with $\tau^*>0$. Thus, we have $z(x,y)>0$ in $\mathbb R^2$. But, that contradicts with $z(0,0)=0$. So, case 2 is ruled out too.

Finally, we have proved that $\tau^*=0$, which means that $\phi\leq\phi^{\tau}$ for all $\tau\geq0$. Then, it follows from the discussion in the beginning of this proof that the function $\phi$ is increasing in $y$. Thus, the proof of Theorem \ref{phi is y increasing} is complete.\hfill$\Box$


\section{Proof of the asymptotic behaviors}\label{proofs of asymptotics}

This section is devoted to the proof of Theorem \ref{c3 asymptotics}. We begin first with Parts~i) and~iii). It follows from formula (\ref{c^*}) that for all $\gamma\geq0$, $m>0$, $\rho>0$ and $L>0$
\begin{equation}\label{parametric M conical speed}\left\{\baa{rcl}
\displaystyle{\frac{c^*(\rho, m^\gamma q,mf)}{\sqrt{m}}} & = & \max\left(\frac{\displaystyle{c^*_{\rho A,m^\gamma q\sin\alpha,mf}}}{\sqrt{m}\sin\alpha},\frac{\displaystyle{c^*_{\rho B,m^\gamma q\sin\beta,mf}}}{\sqrt{m}\sin\beta}\right)\!,\vspace{3pt}\\
\displaystyle{\frac{c^*(m\rho, m^\gamma q,f)}{\sqrt{m}}} & = & \max\left(\frac{\displaystyle{c^*_{m\rho A,m^\gamma q\sin\alpha,f}}}{\sqrt{m}\sin\alpha},\frac{\displaystyle{c^*_{m\rho B,m^\gamma q\sin\beta,f}}}{\sqrt{m}\sin\beta}\right)\!,\vspace{3pt}\\
c^*(\rho,q_L,f) & = & \max\left(\frac{\displaystyle{c^*_{\rho A,q_L\sin\alpha,f}}}{\sin\alpha},\frac{\displaystyle{c^*_{\rho B,q_L\sin\beta,f}}}{\sin\beta}\right)\!.\eaa\right.
\end{equation}
We recall that the quantities appearing in the right-hand side of~(\ref{parametric M conical speed}) are the parametric minimal speeds of propagation of some associated ``left" and ``right" reaction-advection-diffusion problems of the type~(\ref{subbht}). Since $\nabla\cdot Ae=\nabla\cdot Be=0,$ with $e=(0,1)$, in $\mathbb{R}^2,$ $e\cdot \rho Ae=e\cdot\rho Be=\rho$ and the function $f$ satisfies the KPP condition~(\ref{cf}), it follows then from Theorems 4.1, 4.3 and~5.2 of El~Smaily~\cite{El Smaily} that
$$\forall\,\gamma\in[0,1/2],~\lim_{m\rightarrow+\infty}\frac{\displaystyle{c^*_{m\rho A,m^\gamma q\sin\alpha,f}}}{\sqrt{m}}=\lim_{m\rightarrow\infty}\frac{\displaystyle{c^*_{m\rho B,m^\gamma q\sin\beta,f}}}{\sqrt{m}}=2\sqrt{\rho f'(0)},$$
$$\forall\,\gamma\geq1/2,~\lim_{m\rightarrow0^+}\frac{\displaystyle{c^*_{\rho A,m^\gamma q\sin\alpha,\,mf}}}{\sqrt{m}}=\lim_{m\rightarrow0^+}\frac{\displaystyle{c^*_{\rho B,m^\gamma q\sin\beta,\,mf}}}{\sqrt{m}}=2\sqrt{\rho f'(0)},$$
and
$$\lim_{L\rightarrow0^+}\displaystyle{c^*_{\rho A, q_L\sin\alpha,f}}=\lim_{L\rightarrow0^+}\displaystyle{c^*_{\rho B, q_L\sin\beta,f}}=2\sqrt{\rho f'(0)}.$$
Together with (\ref{parametric M conical speed}), we obtain the limits (\ref{limit with small reaction}), (\ref{limit with large diffusion}) and (\ref{limit as L goes to zero}).

Let us now turn to the proof of Part ii) of Theorem \ref{c3 asymptotics}. Remember first that
\be\label{formulem}
\frac{c^*(\rho,mq,f)}{m}=\max\left(\frac{c^*_{\rho A,mq\sin\alpha,f}}{m\sin\alpha},\frac{c^*_{\rho B,mq\sin\beta,f}}{m\sin\beta}\right)
\ee
for all $m>0$, from Theorem~\ref{c3theorem1}. Let now $\widetilde{q}$ be the vector field defined by
$$\widetilde{q}(x,y)=(0,q(x))\ \hbox{ for all }(x,y)\in\mathbb{R}^2.$$
This field is $(L,l)$-periodic in $\mathbb{R}^2$ for each $l>0$, and it satisfies $\nabla\cdot\widetilde{q}=0$ in $\mathbb{R}^2$. Therefore, it follows from Theorem 1.1 in \cite{EK1} or Theorem 1.1 in \cite{zlatos} that, for each $l>0$,
\be\label{limitsl}\left\{\baa{rcl}
\displaystyle{\frac{c^*_{\rho A,mq\sin\alpha,f}}{m\sin\alpha}} & \displaystyle{\mathop{\longrightarrow}_{m\to+\infty}} & \lambda_{\rho A,l},\vspace{3pt}\\
\displaystyle{\frac{c^*_{\rho B,mq\sin\beta,f}}{m\sin\beta}} & \displaystyle{\mathop{\longrightarrow}_{m\to+\infty}} & \lambda_{\rho B,l},\eaa\right.
\ee
where, for any matrix~$M$ fulfilling~(\ref{Mdef+}) and for any $l>0$, the quantity~$\lambda_{M,l}$ is defined by
$$\lambda_{M,l}=\max_{w\in\mathcal{I}_{M,l}}R_{M,l}(w),\qquad R_{M,l}(w)=\frac{\displaystyle{\int_{(0,L)\times(0,l)}}q\,w^2}{\displaystyle{\int_{(0,L)\times(0,l)}}w^2}$$
and
$$\baa{rcl}
\mathcal{I}_{M,l} & = & \Big\{w\in H^1_{loc}(\mathbb{R}^2)\backslash\{0\},\ w\hbox{ is }(L,l)\hbox{-periodic},\ \widetilde{q}\cdot\nabla w=0\hbox{ a.e. in }\mathbb{R}^2,\\
& &  \quad\int_{(0,L)\times(0,l)}\nabla w\cdot M\nabla w\le f'(0)\int_{(0,L)\times(0,l)}w^2\Big\}\eaa$$
is a subset of the set of non-trivial $(L,l)$-periodic first integrals of $\widetilde{q}$. Notice that the set~$\mathcal{I}_{M,l}$ contains the non-zero constants, and that the max in the definition of $\lambda_{M,l}$ is reached, see~\cite{EK1,zlatos}. It follows from~(\ref{limitsl}) that the quantities~$\lambda_{\rho A,l}$ and~$\lambda_{\rho B,l}$ do not depend on $l>0$. Furthermore, since $\widetilde{q}(x,y)=(0,q(x))$, there holds
\be\label{ineqlambda}
\lambda_{\rho A,l}\ge\lambda_{\rho A,0}\ \hbox{ and }\ \lambda_{\rho B,l}\ge\lambda_{\rho B,0}\ \hbox{ for all }l>0,
\ee
where, for any matrix~$M$ fulfilling~(\ref{Mdef+}),
$$\lambda_{M,0}=\max_{\substack{w\in H^1_{loc}(\mathbb{R})\backslash\{0\},\,L\hbox{-periodic}\\ M_{1,1}\|w'\|_{L^2(0,L)}^2\le f'(0)\|w\|_{L^2(0,L)}^2}}\ \frac{\displaystyle{\int_0^L}q\,w^2}{\displaystyle{\int_0^L}w^2}.$$
Let us now check that the opposite inequalities $\lambda_{\rho A,l}\le\lambda_{\rho A,0}$ and $\lambda_{\rho B,l}\le\lambda_{\rho B,0}$ also hold. The proof uses elementary arguments, we just sketch it here for the sake of completeness. We do it for $\lambda_{\rho A,l}$, the proof being identical for $\lambda_{\rho B,l}$. Let $\{l_n\}_{n\in\mathbb{N}}$ be the sequence of positive real numbers defined by $l_n=2^{-n}$ for all $n\in\mathbb{N}$, and let $\{w_n\}_{n\in\mathbb{N}}$ be a sequence of maximizers of the functionals $R_{\rho A,l_n}$ in $\mathcal{I}_{\rho A,l_n}$, that is
\be\label{lambdan}
\lambda_{\rho A,l_n}=R_{\rho A,l_n}(w_n)=\frac{\displaystyle{\int_{(0,L)\times(0,l_n)}}q\,w_n^2}{\displaystyle{\int_{(0,L)\times(0,l_n)}}w_n^2}=\frac{\displaystyle{\int_{(0,L)\times(0,1)}}q\,w_n^2}{\displaystyle{\int_{(0,L)\times(0,1)}}w_n^2}
\ee
for all $n\in\mathbb{N}$. Without loss of generality, one can assume that $\|w_n\|_{L^2((0,L)\times(0,1))}=1$ for all $n\in\mathbb{N}$. By definition, one has
$$\baa{rcl}
\rho\displaystyle{\int_{(0,L)\times(0,1)}}\nabla w_n\cdot A\nabla w_n & = & 2^n\rho\displaystyle{\int_{(0,L)\times(0,l_n)}}\nabla w_n\cdot A\nabla w_n\vspace{3pt}\\
& \le & 2^nf'(0)\displaystyle{\int_{(0,L)\times(0,l_n)}}w_n^2\,=\,f'(0)\displaystyle{\int_{(0,L)\times(0,1)}}w_n^2\,=\,f'(0).\eaa$$
By coercivity of the matrix $A$, the sequence $\{w_n\}_{n\in\mathbb{N}}$ is then bounded in $H^1((0,L)\times(0,1))$. There exists then a function $w_{\infty}\in H^1_{loc}(\mathbb{R}^2)$, which is $(L,1)$-periodic, such that, up to extraction of a sequence, $w_n\to w_{\infty}$ as $n\to+\infty$ in $L^2_{loc}(\mathbb{R}^2)$ strongly and in $H^1_{loc}(\mathbb{R}^2)$ weakly. Thus,
$$\rho\!\int_{(0,L)\times(0,1)}\!\!\nabla w_{\infty}\cdot A\nabla w_{\infty}\le\liminf_{n\to+\infty}\,\rho\!\displaystyle{\int_{(0,L)\times(0,1)}}\!\!\nabla w_n\cdot A\nabla w_n\le f'(0)\!\int_{(0,L)\times(0,1)}\!\!w_{\infty}^2=f'(0).$$
It is then classical to see that $w_{\infty}$ does not depend on $y$. Therefore,
$$\lambda_{\rho A,l_n}\mathop{\longrightarrow}_{n\to+\infty}\frac{\displaystyle{\int_0^L}q\,w_{\infty}^2}{\displaystyle{\int_0^L}w_{\infty}^2}\le\lambda_{\rho A,0}$$
from (\ref{lambdan}) and the definition of $\lambda_{\rho A,0}$. Together with (\ref{ineqlambda}) and the fact that the quantities~$\lambda_{\rho A,l}$ do not depend on~$l$, one concludes that $\lambda_{\rho A,l}=\lambda_{\rho A,0}$ for all $l>0$. It follows then from~(\ref{defAB}),~(\ref{formulem}) and~(\ref{limitsl}) that
$$\frac{c^*(\rho,mq,f)}{m}\mathop{\longrightarrow}_{m\to+\infty}\max\big(\lambda_{\rho A,0},\lambda_{\rho B,0}\big)=\max_{\substack{w\in H^1_{loc}(\mathbb{R})\backslash\{0\},\,L\hbox{-periodic}\\ \rho\|w'\|_{L^2(0,L)}^2\le f'(0)\|w\|_{L^2(0,L)}^2}}\ \frac{\displaystyle{\int_0^L}q\,w^2}{\displaystyle{\int_0^L}w^2}.$$
This provides~(\ref{within large advection}).\par
Formula~(\ref{within large advection}), together with~(\ref{cq}), implies that~(\ref{large advection small reaction}) and~(\ref{mqmuf}) hold, as in~\cite{EK1,zlatos}. The proof of Theorem~\ref{c3 asymptotics} is thereby complete.\hfill$\Box$


\end{document}